\documentclass[10pt]{article}
\usepackage{amsfonts}
\usepackage{hyperref}
\usepackage{mathrsfs,enumerate,amssymb,amsmath,color,lineno}
\usepackage{indentfirst}
\usepackage{amsmath}
\usepackage{amssymb}
\usepackage[amsmath,thmmarks]{ntheorem}
\usepackage{cite}
\usepackage{graphicx}
\usepackage{tikz}
\usepackage{caption}
\usepackage{subcaption}
\usepackage{caption}
\newtheorem{definition}{\bf Definition}[section]
\newtheorem{lemma}{\bf Lemma}[section]
\newtheorem{theorem}{\bf Theorem}[section]
\newtheorem{remark}{\bf Remark}[section]
\newtheorem{corollary}{\bf Corollary}[section]
\newtheorem{example}{\bf Example}[section]

\newtheorem{proposition}{\bf Proposition}[section]

\begin{document}
\setcounter{page}{1}

\title{{\textbf{Orders of continuous cancellative triangular subnorms}}\thanks {Supported by
the National Natural Science Foundation of China (No.12471440)}}
\author{Ting Tang\footnote{\emph{E-mail address}: Tangting199099@cwnu.edu.cn}, Xue-ping Wang\footnote{Corresponding author. xpwang1@hotmail.com; fax: +86-28-84761502},\\
\emph{(School of Mathematical Sciences, Sichuan Normal University,}\\
\emph{Chengdu 610066, Sichuan, People's Republic of China)}}

\newcommand{\pp}[2]{\frac{\partial #1}{\partial #2}}
\date{}
\maketitle
\begin{quote}
{\bf Abstract} We focus on the order relations of continuous cancellative t-subnorms. First, we present some necessary and sufficient conditions along with several interesting sufficient criteria for the comparability of continuous cancellative t-subnorms. Then we characterize the growth and boundedness of additive generators, which are used for the comparison of continuous cancellative t-subnorms.

\end{quote}

\section{Introduction}
Triangular norms (t-norms for short), introduced by Schweizer and Sklar \cite{BS1960} to model the triangle inequality in probabilistic (statistical) metric spaces,  play an indispensable role in domains such as probabilistic metric spaces \cite{BS1983}, fuzzy logic \cite{PH1998}, fuzzy control \cite{MS1985}, non-additive measures and integrals \cite{EP1995}. For further details on t-norms and their applications, see \cite{CA2006,EP2000}. Generalizing t-norms, Jenei \cite{SJ2001} introduced triangular subnorms (t-subnorms for short) as binary operations that need not possess a unit element.

\begin{figure}[!h]
    \centering
    \begin{subfigure}[t]{0.9\textwidth}
        \centering
        \includegraphics[width=\textwidth]{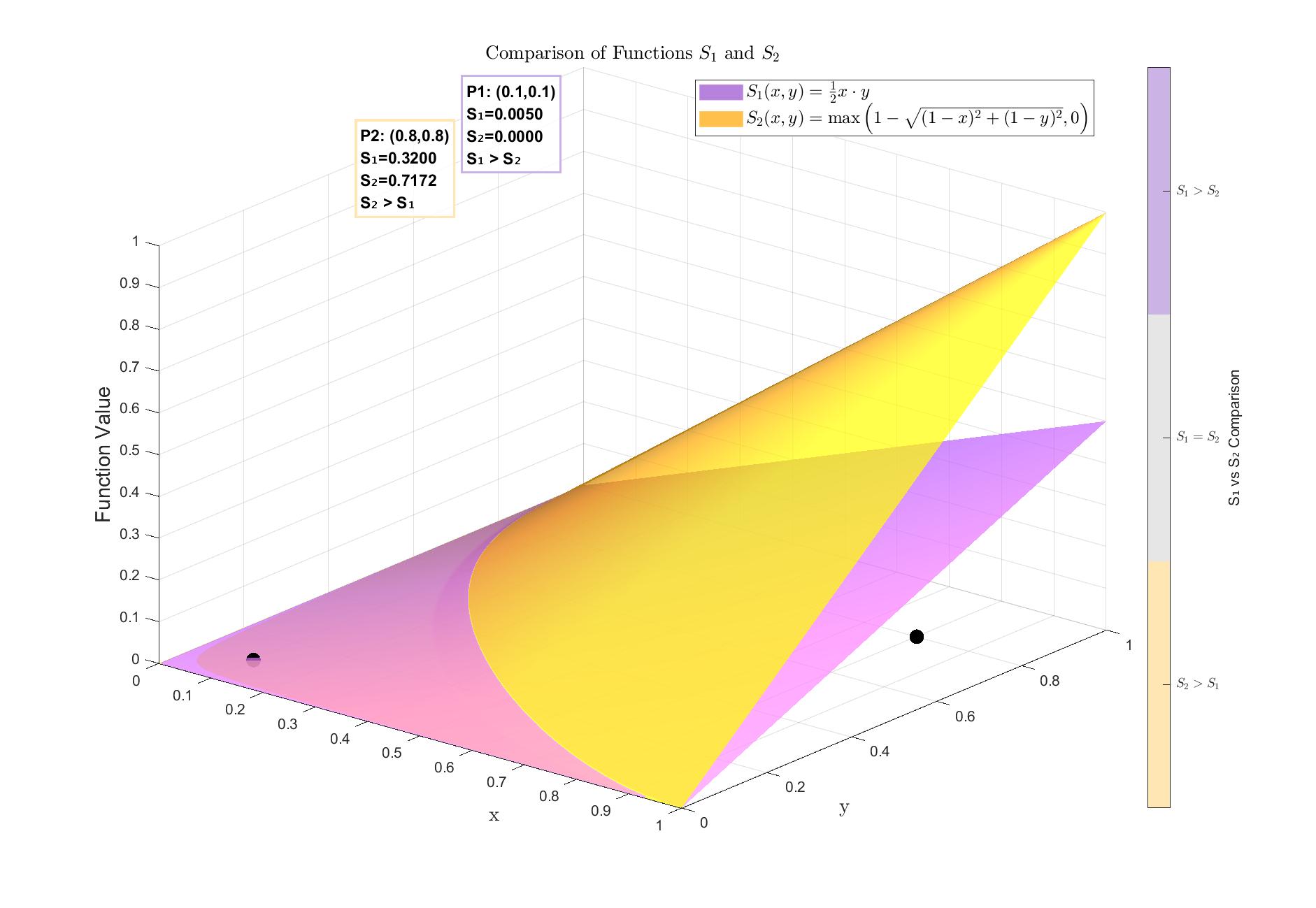}
        \caption{$S_{1}(x,y)=\frac{xy}{2}$ and $S_{2}=T_{2}^{Y}$(Yager t-norm, see \cite{EP2000}).}
    \end{subfigure}
    \hfill
    \begin{subfigure}[t]{0.9\textwidth}
        \centering
        \includegraphics[width=\textwidth]{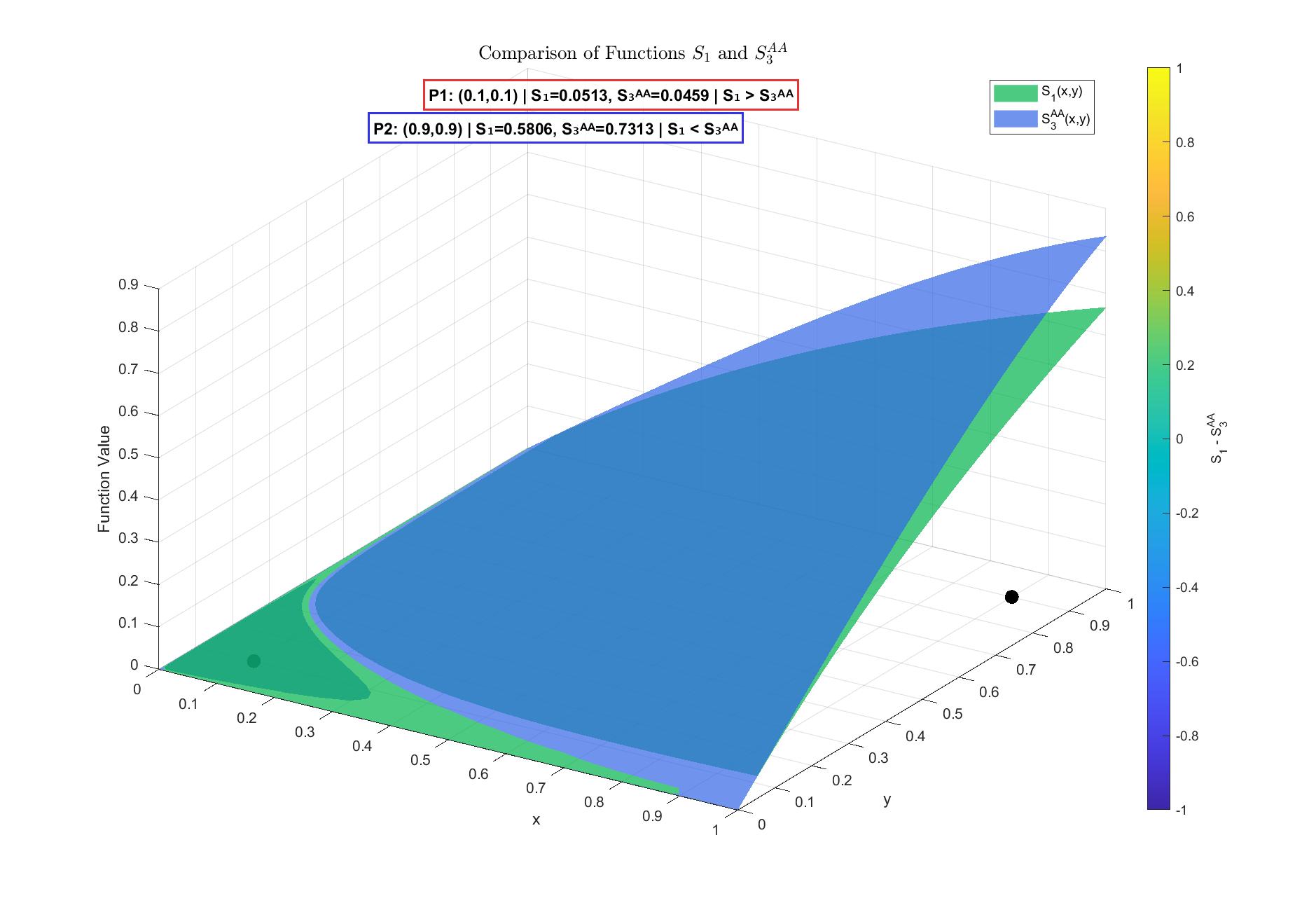}
        \caption{$S_{1}(x,y)=\frac{xy}{x+y-0.5xy}$ and $S_{3}^{AA}$ when $a=0.5$ (see Example \ref{exp3.3} $(\textup{ii})$).}
    \end{subfigure}
    \caption{3D plots of incomparable t-subnorms.}
    \label{fig1}
\end{figure}

As the renowned mathematician Edward V. Huntington stated, ``The fundamental importance of the subject of order may be inferred from the fact that all the concepts required in geometry can be expressed in terms of the concept of order alone". This insight underscores the importance of studying the order relations of t-subnorms. Several contributions to this subject have been made in \cite{BS1983,EP1997}, and these works are about the natural order in the class of continuous (Archimedean) t-norms. However, there is hardly any research directly targeting the order relations of t-subnorms. This leads us to be interesting in the question whether, given two t-subnorms $S_{1}$ and $S_{2}$, $S_{1}$ is weaker than $S_{2}$ or, equivalently, $S_{2}$ is stronger than $S_{1}$ (denoted by $S_{1}\leq S_{2}$), i.e., $S_{1}(x,y)\leq S_{2}(x,y)$ for all $(x,y)\in[0,1]^{2}$. Obviously, the relation $\leq$ is a partial order but not a linear order on the class of all t-subnorms, as there exist incomparable t-subnorms (see Fig. \ref{fig1}). It is well-known that for two given t-subnorms $S_{1}$ and $S_{2}$ it is quite often difficult to directly check whether $S_{1}\leq S_{2}$ or not, even impossible. The main reason stems from the poorly understood structure of t-subnorms as binary functions -- even continuous ones. In fact, for continuous t-subnorms, their structure is only known in highly specific cases: either requiring a neutral element (reducing them to continuous t-norms, see \cite{EP2000}) or cancellativity (enabling the application of Acz\'{e}l's classical results \cite{JA1966}, see \cite{AM2002} for detail). In this article, we focus on comparing continuous cancellative t-subnorms through deeply understanding the structure of them.

The rest of this article is organized as follows. In Section 2, we recall some basic concepts and results of t-norms and t-subnorms, respectively. In Section 3, we explore the comparison of continuous cancellative t-subnorms, establishing some necessary and sufficient conditions along with some distinct sufficient criteria. In Section 4, we characterize the growth and boundedness of additive generators for applying to the comparison of continuous cancellative t-subnorms. A conclusion is drawn in Section 5.

\section{Preliminaries}

In this section, we recall some definitions and results about t-norms and t-subnorms which will be used in the sequel.

\begin{definition}[\cite{EP2000}]\label{def2.1}
\emph{A binary operation $T:[0,1]^{2}\rightarrow[0,1]$ is a $t$-$norm$ if it is commutative, associative, non-decreasing in both variables and 1 is its neutral element.}
\end{definition}

\begin{definition}[\cite{EP2000}]\label{def2.2}
\emph{A binary operation $C:[0,1]^{2}\rightarrow[0,1]$ is a $t$-$conorm$ if it is commutative, associative, non-decreasing in both variables and 0 is its neutral element.}
\end{definition}

The duality between t-norms and t-conorms is expressed by the fact that from any t-norm $T$ we can obtain its dual t-conorm $C$ by the equation $C(x,y)=1-T(1-x,1-y)$ and vice versa.

In literature, the four t-norms: $T_{M}, T_{P}, T_{L}$, and $T_{D}$ are often discussed (see \cite{EP2000}). Note that $T_{D}<T_{L}<T_{P}<T_{M}$ and $T_{D}\leq T\leq T_{M}$  for an arbitrary t-norm $T$.

\begin{definition}[\cite{EP2000}]\label{def:2.3}
 \emph{A t-norm $T:[0,1]^{2}\rightarrow[0,1]$ is said to be
 \begin{enumerate}
\renewcommand{\labelenumi}{(\roman{enumi})}
\item \emph{Cancellative} if $T(x,y)=T(x,z)$ for all $x,y,z\in[0,1]$ with $x>0$ then $y=z$;
  \item \emph{Archimedean} if for every $(x,y)\in(0,1)^2$ there exists an $n\in \mathbb{N}=\{1,2,\cdots, n,\cdots\}$ such that $x_{T}^{(n)}<y$, where $x_{T}^{(n)}=T\underbrace{(x,T(x,\cdots))}_{n\  times}$;
  \item \emph{Strict} if it is continuous and strictly monotone;
  \item \emph{Nilpotent} if it is continuous and for each $x\in(0,1)$ there exists an $n\in\mathbb{N}$ such that $x_{T}^{(n)}=0$.
\end{enumerate}}
\end{definition}

It is well known that a t-norm $T$ is strictly monotone if and only if it satisfies the cancellation law, and each continuous Archimedean t-norm (resp. t-conorm) is either strict or nilpotent.

\begin{definition}[\cite{EP2000,EP1999}]\label{def2.4}
\emph{Let $a, b, c, d\in [-\infty, \infty]$ with $a<b, c<d$ and $f:[a,b]\rightarrow[c,d]$ be a non-increasing function. Then the function $f^{(-1)}:[c,d]\rightarrow[a,b]$ defined by
\begin{equation*}
f^{(-1)}(y)=\sup\{x\in [a,b]\mid f(x)>y\}
\end{equation*}
is called a \emph{pseudo-inverse} of a non-increasing function $f$.}
\end{definition}

The following is the representation theorem of Archimedean t-norms \cite{EP2000,EP1999}.

\begin{theorem}\label{theorem:2.1}
For a function $T:[0,1]^{2}\rightarrow[0,1]$ the following are equivalent:
\begin{enumerate}
\renewcommand{\labelenumi}{(\roman{enumi})}
\item $T$ is a continuous Archimedean t-norm.
\item $T$ has a continuous additive generator, i.e., there exists a continuous, strictly decreasing function $t:[0,1]\rightarrow[0,\infty]$ with $t(1)=0$, which
is uniquely determined up to a positive multiplicative constant, such that for all $(x,y)\in[0,1]^{2}$ there is
\begin{equation*}
T(x, y) = t^{(-1)}(t(x) + t(y))
\end{equation*}
where $t^{(-1)}$ is the pseudo-inverse of t, i.e.,
\begin{equation*}
t^{(-1)}(x) =
\begin{cases}
t^{-1}(x), & 0 \leq x \leq t(0), \\
0, & x \geq t(0).
\end{cases}
\end{equation*}
\end{enumerate}
\end{theorem}

The function $t$ in Theorem \ref{theorem:2.1} is called an \emph{additive generator} of the t-norm $T$. If $t(0)=\infty$ then $t^{(-1)}$ is the inverse of $t$ and $T$ is strict, if $t(0)<\infty$ then $T$ is nilpotent.

Next we recall the definitions of t-subnorm and its dual t-superconorm, respectively.

\begin{definition}[\cite{SJ2001}]\label{def2.5}
 \emph{A binary operation $S:[0, 1]^2\rightarrow[0, 1]$ is a $t$-$subnorm$ if it is commutative, associative, non-decreasing in both variables and $S(x, y)\leq \min(x, y)$ for all $(x, y)\in[0, 1]^2$.}
\end{definition}

\begin{definition}[\cite{AM2016}]\label{def2.6}
\emph{A binary operation $M:[0, 1]^2\rightarrow[0, 1]$ is a $t$-$superconorm$ if it is commutative, associative, non-decreasing in both variables and $M(x, y)\geq \max(x, y)$ for all $(x, y)\in[0, 1]^2$.}
\end{definition}

Between t-subnorms and t-superconorms there is the same duality as between t-norms and t-conorms. Thus all results which we shall show for t-subnorms can be obtained immediately also for t-superconorms.

Evidently, each t-norm (resp. t-conorm) is also a t-subnorm (resp. t-superconorm). When a t-subnorm is not a t-norm, it is called a proper t-subnorm. In particular, Mesiarov\'{a}-Zem\'{a}nkov\'{a} \cite{AM2004} proved that a continuous t-subnorm $S$ is proper if and only if $S(1,1)<1$.

The zero t-subnorm $Z$, where $Z(x, y)=0$ for all $(x,y)\in[0,1]^{2}$, is the weakest t-subnorm. Note that $Z\leq S\leq T_{M}$ for an arbitrary  t-subnorm $S$.

\begin{definition}[\cite{AM2016}]\label{def2.7}
 \emph{A t-subnorm $S:[0, 1]^2\rightarrow[0, 1]$ is said to be
 \begin{enumerate}
\renewcommand{\labelenumi}{(\roman{enumi})}
\item \emph{Cancellative} if $S(x, y)=S(x, z)$ for all $x, y, z\in[0, 1]$ with $x>0$ implies $y=z$;
  \item \emph{Archimedean} if for all $x, y\in(0, 1)$  there exists an $n\in\mathbb{N}$ such that $x_{S}^{(n)}<y$, where $x_{S}^{(n)}=S\underbrace{(x,S(x,\cdots))}_{n\  times}$.
\end{enumerate}}
\end{definition}

Observe that the fact that $1$ is the neutral element of t-norms forces the strict monotonicity of any additive generator $t:[0, 1]\rightarrow[0, \infty]$ of a t-norm as well as $t(1)=0$. However, for t-subnorms this restriction can be relaxed.

\begin{definition}[\cite{AM2004}]\label{def2.8}
\emph{A non-increasing function $s:[0, 1]\rightarrow[0, \infty]$ is called an \emph{additive generator} of a t-subnorm $S:[0, 1]^2\rightarrow[0, 1]$ if for all $(x, y)\in[0, 1]^2$,
\begin{equation*}
S(x, y)=s^{(-1)}(s(x) + s(y))
\end{equation*} where $s^{(-1)}$ is the pseudo-inverse of s.}
\end{definition}

We only mention that, the continuity of an additive generator $s:[0, 1]\rightarrow[0, \infty]$ of a t-subnorm $S$ ensures neither the continuity nor Archimedean property of $S$ (see Proposition 1 and Theorem 3 of \cite{AM2004}). In particular, Mesiarov\'{a} proved the result `` let $S:[0, 1]^2\rightarrow[0, 1]$ be a continuous, Archimedean, proper t-subnorm which is cancellative on $(0, 1]^2$. Then $S$ has a continuous additive generator" (see Theorem 27 of \cite{AM2002}). In this article, we focus on continuous cancellative t-subnorms, which admit a representation theorem stated as follows.
\begin{theorem}[\cite{AM2016}]\label{theorem:2.3}
A t-subnorm $S:[0, 1]^2\rightarrow[0, 1]$ has a continuous, strictly decreasing additive generator $s:[0, 1]\rightarrow[0, \infty]$ with $s(0)=\infty $ if and only if it is continuous and cancellative.
\end{theorem}

Moreover, the additive generator $s$ in Theorem \ref{theorem:2.3} is unique up to a positive multiplicative constant. If $s(1)=0$, then $s:[0,1]\rightarrow[0,\infty]$ is a strictly decreasing bijection and $S$ is a strict t-norm. If $s(1)>0$, then $S$ is a proper t-subnorm, and in this case the pseudo-inverse $s^{(-1)}$ of $s$ is
\begin{equation*}
s^{(-1)}(u)=
\begin{cases}
1, &  0 \leq u \leq s(1), \\
s^{-1}(u), & s(1) \leq u \leq \infty.
\end{cases}
\end{equation*}
Therefore, each continuous cancellative t-subnorm $S$ can be represented in the form
\begin{equation}
\label{eq:1.1}
  S(x, y)=s^{-1}(s(x) + s(y))
\end{equation}
where $s:[0,1]\rightarrow[0,\infty]$ is strictly decreasing continuous with $s(0)=\infty$.

\section{Comparison of continuous cancellative t-subnorms based on the subadditivity}
In this section, we explore the comparison of continuous cancellative t-subnorms through the subadditivity of additive generators.

Recall that additive generators of continuous cancellative t-subnorms contain all the information because of Theorem \ref{theorem:2.3}. We naturally guess that the pointwise comparability of continuous cancellative t-subnorms may also be deducible from the corresponding properties of their respective additive generators. Observe that each additive generator $s$ of a continuous cancellative t-subnorm, viewed as a function from $[0,1]$ to $[s(1),\infty]$, is a strictly decreasing bijection, and subsequently, its inverse function $s^{-1}: [s(1),\infty]\rightarrow[0,1]$ always exists. Write
$$\mathcal{S}=\{S| S:[0, 1]^2\rightarrow[0, 1] \mbox{ is a continuous cancellative t-subnorm}\},$$
$$\mathcal{S}_{p}=\{S| S:[0, 1]^2\rightarrow[0, 1] \mbox{ is a continuous cancellative proper t-subnorm}\} \mbox{ and}$$
$$\mathcal{T}_{s}=\{S| S:[0, 1]^2\rightarrow[0, 1] \mbox{ is a continuous cancellative t-norm}\}.$$
Obviously, both $\mathcal{S}_{p}$ and $\mathcal{T}_{s}$ are proper subsets of $\mathcal{S}$, and their union is $\mathcal{S}$, i.e., $\mathcal{S}=\mathcal{S}_{p}\bigcup\mathcal{T}_{s}$. First, we have the following theorem.
\begin{theorem}\label{theorem:3.1}
Let $S_{1}, S_{2}\in \mathcal{S}$ with additive generators $s_{1}:[0,1]\rightarrow[s_{1}(1), \infty]$ and $s_{2}:[0,1]\rightarrow[s_{2}(1), \infty]$, respectively. Then $S_{1}\leq S_{2}$ if and only if the function $s_{1}\circ s_{2}^{-1}:[s_{2}(1), \infty]\rightarrow[s_{1}(1), \infty]$ is subadditive, i.e.,
\begin{equation*}
s_{1}\circ s_{2}^{-1}(u+v)\leq s_{1}\circ s_{2}^{-1}(u)+s_{1}\circ s_{2}^{-1}(v)
\end{equation*}for all $u,v\in[s_{2}(1),\infty]$.
\end{theorem}
\begin{proof}
Note that the function $s_{1}\circ s_{2}^{-1}:[s_{2}(1), \infty]\rightarrow[s_{1}(1), \infty]$ is continuous and strictly increasing bijection since $S_{1}, S_{2}\in \mathcal{S}$. By Eq.(\ref{eq:1.1}), $S_{1}\leq S_{2}$ if and only if
\begin{equation} \label{eq:1.2}
s_1^{-1}(s_1(x) + s_1(y))\leq s_2^{-1}(s_2(x) + s_2(y))
\end{equation} for all $x,y\in[0,1]$.
Fix $x,y\in[0,1]$, and set $u=s_{2}(x)$ and $v=s_{2}(y)$. Taking into account that $s_{1}$ is strictly decreasing and applying $s_{1}$ to
both sides of Eq.(\ref{eq:1.2}), we have that Eq.(\ref{eq:1.2})is equivalent to the fact that
$$s_{1}\circ s_{2}^{-1}(u)+s_{1}\circ s_{2}^{-1}(v)\geq s_{1}\circ s_{2}^{-1}(u+v)$$
 for all $u,v\in[s_{2}(1),\infty]$. Therefore, letting $h=s_{1} \circ s_{2}^{-1}$, we get $S_{1}\leq S_{2}$ if and only if $$h(u+v)\leq h(u)+h(v)$$
 for all $u,v\in[s_{2}(1),\infty]$, i.e., $h$ is subadditive.
\end{proof}

For any continuous cancellative proper t-subnorm $S$, the additive generator $s$ associated with $S$ satisfies the condition $s(1)>0$. Thus for the additive generator $s$ of any continuous cancellative proper t-subnorm $S$, the function $\frac{s}{s(1)}:[0,1]\rightarrow[1,\infty]$ is uniquely determined and it is called a \emph{normalized additive generator} of $S$.

In what follows, without loss of generality, every additive generator of $S$ always refers to its normalized additive generator when $S\in\mathcal{S}_{p}$.

From Theorem \ref{theorem:3.1}, the following remark holds.
\begin{remark}\label{remark:3.1} \emph{Let $S_{1}, S_{2}\in\mathcal{S}$ with additive generators  $s_{1}:[0,1]\rightarrow[s_{1}(1), \infty]$ and $s_{2}: [0,1]\rightarrow[s_{2}(1), \infty]$, respectively.
\renewcommand{\labelenumi}{(\roman{enumi})}
\begin{enumerate}
\item If $s_{1}(1)=s_{2}(1)=0$, then $S_{1}, S_{2}\in\mathcal{T}_{s}$, and $S_{1}\leq S_{2}$ if and only if the function $s_{1}\circ s_{2}^{-1}: [0,\infty]\rightarrow[0,\infty]$ is subadditive. Thus Theorem \ref{theorem:3.1} generalizes Theorem 7 of \cite{BS1961}.
\item If $s_{1}(1)=1$ and $s_{2}(1)=0$, then $S_{1}\in\mathcal{S}_{p}$, $S_{2}\in\mathcal{T}_{s}$, and $S_{1}\leq S_{2}$ if and only if the function $s_{1}\circ s_{2}^{-1}:[0,\infty]\rightarrow[1,\infty]$ is subadditive.
\item If $s_{1}(1)=s_{2}(1)=1$, then $S_{1}, S_{2}\in\mathcal{S}_{p}$, and $S_{1}\leq S_{2}$ if and only if the function $s_{1}\circ s_{2}^{-1}: [1,\infty]\rightarrow[1,\infty]$ is subadditive.
\item There are no $S\in\mathcal{S}_{p}$ and $T\in\mathcal{T}_{s}$ satisfying $T\leq S$. In fact, for every t-subnorm $S$ and t-norm $T$, we have $S(x,1)\leq \min (x,1)=T(x,1)$ for all $x\in[0,1]$.
\end{enumerate}}
\end{remark}

The following example illustrates Theorem \ref{theorem:3.1}.

\begin{example}\label{exp3.1}\emph{
Considering the Hamacher product, i.e., $T_{0}^{H}$ (taking into account $\frac{1}{\infty}=0$ and $\frac{1}{0}=\infty$), we get
\begin{equation*}
T_{0}^{H}(x,y)=\frac{1}{\frac{1}{x}+\frac{1}{y}-1}
\end{equation*} for all $(x,y)\in[0,1]^{2}$.
Moreover, the additive generator $t:[0, 1]\rightarrow[0, \infty] $ of $T_{0}^{H}$ is $t(x)=\frac{1-x}{x}$.
\renewcommand{\labelenumi}{(\roman{enumi})}
\begin{enumerate}
\item Let $T_{1}=T_{p}$ and $T_{2}=T_{0}^{H}$. Then their additive generators are $t_{1}(x)=-\ln x$ and $t_{2}(x)=\frac{1}{x}-1$, respectively. Clearly, $t_{1}(1)=t_{2}(1)=0$ and $t_{1}\circ t_{2}^{-1}(x)=\ln(x+1)$ for all $x\in[0,\infty]$. Thus the function $t_{1}\circ t_{2}^{-1}: [0, \infty]\rightarrow[0, \infty]$ is subadditive. Therefore, by Remark \ref{remark:3.1}(i) we have $T_{1}\leq T_{2}$.
\item Let $S:[0, 1]^2\rightarrow[0, 1]$ with $S(x,y)=\frac{1}{2}xy$ and $T =T_{0}^{H}$. Then the normalized additive generator of $S$ is the function $s:[0, 1] \rightarrow[1, \infty]$ with $s(x)=-\frac{\ln x}{\ln 2}+1$. Obviously, $s(1)=1$, $t(1)=0$ and $s\circ t^{-1}(x)=\frac{\ln(x+1)}{\ln2}+1$ for all $x\in[0,\infty]$. Thus the function $s\circ t^{-1}:[0, \infty]\rightarrow[1, \infty]$ is subadditive. Therefore, by Remark \ref{remark:3.1}(ii) we have $S\leq T$.
\item Let $S_1: [0, 1]^2\rightarrow[0, 1]$ with $S_{1}(x,y)=\frac{xy}{x+y-0.5xy}$ and $S_2: [0, 1]^2\rightarrow[0, 1]$ with $S_{2}(x,y)=\frac{xy}{x+y-0.7xy}$. Then their normalized additive generators are $s_1:[0, 1]\rightarrow[1, \infty]$ with $s_{1}(x)=\frac{2}{x}-1$ and $s_2:[0, 1]\rightarrow[1, \infty]$ with $s_{2}(x)= \frac{10}{3x}-\frac{7}{3}$, respectively. Evidently, $s_1(1)=s_2(1)=1$ and $s_{1}\circ s_{2}^{-1}(x)=\frac{3x+2}{5}$ for all $x\in[1,\infty]$. Thus the function $s_{1}\circ s_{2}^{-1}:[1, \infty]\rightarrow[1, \infty]$ is subadditive. Therefore, by Remark \ref{remark:3.1}(iii) we have $S_{1}\leq S_{2}$.
\end{enumerate}}
\end{example}

\begin{theorem}\label{theorem:3.2}
Let $S_{1}, S_{2}\in\mathcal{S}$ with additive generators $s_{1}:[0,1]\rightarrow[s_{1}(1), \infty]$ and $s_{2}:[0,1]\rightarrow[s_{2}(1), \infty]$, respectively. Then $S_{1}=S_{2}$ if and only if the function $s_{1}\circ s_{2}^{-1}:[s_{2}(1), \infty]\rightarrow[s_{1}(1), \infty]$ is homogeneous linear, i.e., $s_{1} \circ s_{2}^{-1}(x)=cx$ for all $x\in[s_{2}(1), \infty]$ where $c$ is an arbitrary constant with $c\in(0,\infty)$.
\end{theorem}
\begin{proof}
By Eq.(\ref{eq:1.1}), $S_{1}=S_{2}$ if and only if for all $x,y\in[0,1]$
\begin{equation} \label{eq:1.3}
s_1^{-1}(s_1(x) + s_1(y)) = s_2^{-1}(s_2(x) + s_2(y)).
\end{equation}
Fixing $x,y\in[0,1]$, setting $u=s_{2}(x)$ and $v=s_{2}(y)$, and applying $s_{1}$ to
both sides of Eq.(\ref{eq:1.3}), we have that Eq.(\ref{eq:1.3}) is equivalent to
\begin{equation} \label{eq:1.4}
s_{1}\circ s_{2}^{-1}(u)+s_{1}\circ s_{2}^{-1}(v)=s_{1}\circ s_{2}^{-1}(u+v)
\end{equation}
 for all $u,v\in[s_{2}(1),\infty]$.

 Now, suppose that $S_{1}=S_{2}$. Then Eq.(\ref{eq:1.4}) holds and it is a Cauchy functional equation, whose continuous, strictly increasing solutions $s_1\circ s_2^{-1}:[s_{2}(1), \infty]\rightarrow[s_{1}(1), \infty]$ must satisfy $s_1\circ s_2^{-1}=c\cdot \mathrm{id}_{[s_{2}(1), \infty]}$ for some $c\in(0, \infty)$. As a consequence, $s_{1} \circ s_{2}^{-1}(x)=c\cdot x$ for all $x\in[s_{2}(1), \infty]$ and for some $c\in(0, \infty)$. Therefore, $s_{1} \circ s_{2}^{-1}$ is homogeneous linear.

 Conversely, assume that $s_{1} \circ s_{2}^{-1}(x)=cx$ for all $x\in[s_{2}(1), \infty]$ where $c$ is an arbitrary constant with $c\in(0,\infty)$. Then from Eq.(\ref{eq:1.4}) we get $S_{1}=S_{2}$.
\end{proof}

From Theorem \ref{theorem:3.2}, we have the following remark.

\begin{remark}\label{remark:3.2} \emph{Let $S_{1}, S_{2}\in\mathcal{S}$ with additive generators $s_{1}:[0,1]\rightarrow[s_{1}(1), \infty]$ and $s_{2}: [0,1]\rightarrow[s_{2}(1), \infty]$, respectively.
\renewcommand{\labelenumi}{(\roman{enumi})}
\begin{enumerate}
\item If $S_{1}, S_{2}\in$ $\mathcal{T}_{s}$, then $S_{1}=S_{2}$ if and only if the function $s_{1}\circ s_{2}^{-1}:[0,\infty]\rightarrow[0,\infty]$ is homogeneous linear.
\item If $S_{1}, S_{2}\in\mathcal{S}_{p}$, then $S_{1}=S_{2}$ if and only if the function $s_{1}\circ s_{2}^{-1}:[1,\infty]\rightarrow[1,\infty]$ is homogeneous linear.
\end{enumerate}}
\end{remark}

Theorem \ref{theorem:3.2} indicates that for any $S_{1}, S_{2}\in\mathcal{S}$ if their respective additive generators $s_{1}:[0,1]\rightarrow[s_{1}(1), \infty]$ and $s_{2}: [0,1]\rightarrow[s_{2}(1), \infty]$ satisfy $s_{1}=c\cdot s_{2}$ where $c$ is an arbitrary constant with $c\in(0,\infty)$ then $S_{1}=S_{2}$. Next, we consider the case of non-homogeneous linear relation $s_{1}=c\cdot s_{2}+b$ where $c$ is an arbitrary constant with $c\in(0,\infty)$ and $b\neq 0$. We have the following two propositions.
\begin{proposition}\label{proo3.1}
Let $S_{1}, S_{2}\in\mathcal{S}_{p}$ with additive generators $s_{1}, s_{2}: [0,1]\rightarrow[1, \infty]$, respectively. If $s_{1}=c\cdot s_{2}+b$ where $0<c\leq 1$ and $c+b=1$, then $S_{1}\leq S_{2}$.
\end{proposition}
\begin{proof}
Since $s_1, s_2:[0,1]\rightarrow[1,\infty]$ are continuous, strictly decreasing and bijective, by setting $x=1$, $s_1=c\cdot s_2 + b$ yields $c+b=1$. Thus $s_{1}=c\cdot s_{2}+1-c$ and for all $x,y\in [0,1]$,
\begin{align*}
S_{1}(x,y)&=s_{1}^{-1}(s_{1}(x)+ s_{1}(y))\\
&=s_{1}^{-1}(c\cdot s_{2}(x)+1-c+c\cdot s_{2}(y)+1-c)\\
&=s_{2}^{-1}(s_{2}(x)+s_{2}(y)+\frac{1}{c}-1)\\
&\leq s_{2}^{-1}(s_{2}(x)+s_{2}(y))\mbox{ since }0<c\leq 1\mbox{ and }s_{2}^{-1} \mbox{ is strictly decreasing}\\
&=S_{2}(x,y).
\end{align*}
Therefore, $S_{1}\leq S_{2}$.
\end{proof}

\begin{proposition}\label{proo3.2}
Let $S\in\mathcal{S}_{p}$, $T\in\mathcal{T}_{s}$ with additive generators $s:[0,1]\rightarrow[1, \infty]$ and $t:[0,1]\rightarrow[0, \infty]$, respectively. If $s=c\cdot t+1$ where $c\in(0,\infty)$, then $S\leq T$.
\end{proposition}
\begin{proof}
Because of $s=c\cdot t+1$, we have that for all $x,y\in [0,1]$,
\begin{align*}
S(x,y)&=s^{-1}(s(x)+ s(y))\\
&=s^{-1}(c\cdot t(x)+1+c\cdot t(y)+1)\\
&=t^{-1}(t(x)+t(y)+\frac{1}{c})\\
&\leq t^{-1}(t(x)+t(y))\mbox{ since }0<c< \infty\mbox{ and }t^{-1} \mbox{ is strictly decreasing}\\
&=T(x,y),
\end{align*}
which follows that $S\leq T$.
\end{proof}

Both Propositions \ref{proo3.1} and \ref{proo3.2} are crucial in the construction of two classes of continuous cancellative proper t-subnorms: one consists of same-type operators which are weaker than a given continuous cancellative proper t-subnorm, and the other is of those which are weaker than a given strict t-norm.

 Because the subadditivity of $s_{1}\circ s_{2}^{-1}$ is very difficult to check directly, we try to seek another method for comparing two t-subnorms although Theorem \ref{theorem:3.1} is applied to comparing t-subnorms via the subadditivity of the function $s_{1}\circ s_{2}^{-1}$. First, it is important to note that we cannot directly conclude $S_{1}\leq S_{2}$ just from $s_{1}\leq s_{2}$ even if their respective additive generators $s_{1}$ and $s_{2}$ contain all their information. As shown in subsequent examples, even for two strict t-norms $T_{1}$ and $T_{2}$, the inequality $t_{1}\leq t_{2}$ where $t_{1}$ and $t_{2}$ are their respective additive generators does not imply $T_{1}\leq T_{2}$, and vice versa. In general, the order of two additive generators has no connection with that of their induced t-subnorms.

\begin{example}\label{exp3.2}\emph{
\renewcommand{\labelenumi}{(\roman{enumi})}
\begin{enumerate}
\item Let $t_1:[0, 1]\rightarrow[0, \infty]$ with $t_{1}(x)=\frac{1}{x}-1$ and $t_2:[0, 1]\rightarrow[0, \infty]$ with $t_{2}(x)=\frac{1}{x}-x$. Obviously, $t_{1}\leq t_{2}$. From Theorem \ref{theorem:2.1}, we have that for all $x,y\in [0,1]$,
$$T_{1}(x,y)=t_{1}^{-1}(t_{1}(x)+t_{1}(y))=\frac{xy}{x+y-xy} \mbox{ and}$$
$$T_{2}(x,y)=t_{2}^{-1}(t_{2}(x)+t_{2}(y))=\frac{\sqrt{(\frac{1}{x}-x+\frac{1}{y}-y)^{2}+4}-(\frac{1}{x}-x+\frac{1}{y}-y)}{2}.$$
By leting $x=y=0.5$, we easily verify that
$$T_{1}(0.5,0.5)=\frac{1}{3}>\frac{\sqrt{13}-3}{2}=T_{2}(0.5,0.5),$$
showing that $T_{1}\nleq T_{2}$.
\item Let $T_1=T_{p}$ and $T_{2}(x,y)=T_{2}^{AA}(x,y)=e^{-\sqrt{(-\ln x)^{2}+(-\ln y)^{2}}}$ (i.e., Acz\'{e}l-Alsina t-norms with parameter $\lambda=2$, see \cite{EP2000}). A simple calculation shows that $T_{1}\leq T_{2}$ and $t_{1}\nleq t_{2}$ where $t_1:[0, 1]\rightarrow[0, \infty]$ with $t_{1}(x)=-\ln x$ and $t_2:[0, 1]\rightarrow[0, \infty]$ with $t_{2}(x)=(-\ln x)^{2}$ are their respective additive generators. In fact, $t_{1}(x)\geq t_{2}(x)$ for all $x\in[e^{-1},1]$.
\end{enumerate}}
\end{example}

 Fortunately, we have several criteria used for the comparison of continuous cancellative t-subnorms which are easier to check than the subadditivity in Theorem \ref{theorem:3.1} sometimes. Let's first recall that a function $f:I\rightarrow\mathbb{R}$, where $I\subseteq\mathbb{R}$ is an arbitrary non-empty interval, is called a \emph{concave function} if for all $x,y\in I$ and $\lambda\in[0,1]$ we have $f(\lambda x+(1-\lambda)y)\geq \lambda f(x)+(1-\lambda)f(y)$. Then we have the following proposition.

\begin{proposition}\label{prop3.1}
Let $S_{1}, S_{2}\in \mathcal{S}$ with additive generators $s_{1}:[0,1]\rightarrow[s_{1}(1), \infty]$ and $s_{2}:[0,1]\rightarrow[s_{2}(1), \infty]$, respectively. If the function $s_{1}\circ s_{2}^{-1}:[s_{2}(1), \infty]\rightarrow[s_{1}(1), \infty] $ is concave and satisfies that $s_{1}\circ s_{2}^{-1}(u)\leq u$ holds for all $u\in[1,\infty]$ when both $s_{1}$ and $s_{2}$ are normalized additive generators, then $S_{1}\leq S_{2}$.
\end{proposition}
\begin{proof}
Letting $h=s_{1}\circ s_{2}^{-1}$, we distinguish three cases as below.
\renewcommand{\labelenumi}{(\roman{enumi})}
\begin{enumerate}
\item Let $s_{1}(1)=s_{2}(1)=0$. Then $S_{1}, S_{2}\in$ $\mathcal{T}_{s}$ and $h(0)=s_{1}\circ s_{2}^{-1}(0)=0$.
\begin{itemize}
  \item If $u=v=0$, then it is trivial that $h(u+v)\leq h(u)+h(v)$.
  \item If $u,v\in[0,\infty]$ with $0<u+v<\infty$, then let $\lambda=\frac{u}{u+v}$. Thus we have
  $$ h(u)=h(\frac{u}{u+v}\cdot(u+v)+\frac{v}{u+v}\cdot 0)\geq \frac{u}{u+v}\cdot h(u+v)+\frac{v}{u+v}\cdot h(0)$$
  since $h$ is a concave function, furthermore, it follows from $h(0)=0$ that
  \begin{equation} \label{eq:1.6}
  h(u)\geq\frac{u}{u+v}\cdot h(u+v).
  \end{equation}
  Similarly, by letting $\lambda=\frac{v}{u+v}$ we have
  \begin{equation*}
  h(v)\geq\frac{v}{u+v}\cdot h(u+v),
  \end{equation*}
 which together with Eq.(\ref{eq:1.6}) yields $h(u+v)\leq h(u)+h(v)$.
  \item If $u+v=\infty$, then either $u$ or $v$ are infinite, and as a consequence,
  $$h(u+v)=h(\infty) \leq h(\infty)+h(\min(u, v))=h(u)+h(v).$$
\end{itemize}
\item Let $s_{1}(1)=1$ and $s_{2}(1)=0$. Then $S_{1}\in\mathcal{S}_{p}$, $S_{2}\in\mathcal{T}_{s}$ and $h(0)=s_{1}\circ s_{2}^{-1}(0)=1$. Set $\lambda=\frac{u}{u+v}$ and $\lambda=\frac{v}{u+v}$, respectively. Then we similarly get $$ h(u)\geq \frac{u}{u+v}\cdot h(u+v)+\frac{v}{u+v}\cdot h(0)=\frac{u}{u+v}\cdot h(u+v)+\frac{v}{u+v}\geq \frac{u}{u+v}\cdot h(u+v)$$ and
   $$ h(v)\geq \frac{v}{u+v}\cdot h(u+v)+\frac{u}{u+v}\cdot h(0)=\frac{v}{u+v}\cdot h(u+v)+\frac{u}{u+v}\geq \frac{v}{u+v}\cdot h(u+v),$$respectively. The rest of the proof is analogous to $(\textup{i})$.
\item Let $s_{1}(1)=s_{2}(1)=1$. Then $S_{1}, S_{2}\in$ $\mathcal{S}_{p}$ and $h(1)=s_{1}\circ s_{2}^{-1}(1)=1$ owing to the fact that both $s_{1}$ and $s_{2}$ are normalized additive generators.
\begin{itemize}
  \item If $u=v=1$, then it is clear that $h(u+v)=h(2)\leq 2=h(u)+h(v)$ since $h(u)\leq u$ holds for all $u\in[1,\infty]$.
  \item If $u,v\in[1,\infty]$ with $1<u+v<\infty$, then let $\lambda=\frac{u-1}{u+v-1}$. Hence
  $$ h(u)=h(\frac{u-1}{u+v-1}\cdot(u+v)+\frac{v}{u+v-1}\cdot 1)\geq \frac{u-1}{u+v-1}\cdot h(u+v)+\frac{v}{u+v-1}\cdot h(1)$$
  since $h$ is a concave function.
  Thus
  \begin{equation*}
  h(u)\geq\frac{u-1}{u+v-1}\cdot h(u+v)+\frac{v}{u+v-1}
  \end{equation*} since $h(1)=1$.
  Analogously, letting $\lambda=\frac{v-1}{u+v-1}$, we have
  \begin{equation*}
  h(v)\geq\frac{v-1}{u+v-1}\cdot h(u+v)+\frac{u}{u+v-1},
  \end{equation*}
  and consequently,
  $$ h(u)+h(v)\geq \frac{u+v-2}{u+v-1}\cdot h(u+v)+\frac{u+v}{u+v-1}.$$
  Then
  $$ h(u)+h(v)-h(u+v)\geq \frac{-1}{u+v-1}\cdot h(u+v)+\frac{u+v}{u+v-1}= \frac{(u+v)-h(u+v)}{u+v-1}.$$
  Thus $h(u)+h(v)-h(u+v)\geq 0$ since $h(u)\leq u$ holds for all $u\in[1,\infty]$. Therefore, $h(u+v)\leq h(u)+h(v)$.
  \item If $u+v=\infty$, then similarly to $(\textup{i})$, we have $h(u+v)\leq h(u)+h(v)$.
\end{itemize}
\end{enumerate}
(i), (ii) and (iii) imply that $h$ is a subadditive function, and thus $S_{1}\leq S_{2}$ from Theorem \ref{theorem:3.1}.
\end{proof}

Generally, the converse of Proposition \ref{prop3.1} does not hold because the subadditivity of a function does not imply its concavity. For example, let $h:[0,\infty]\rightarrow[0,\infty]$ be given by
\begin{equation*}
 h(x)=\begin{cases}
2x & \hbox{if }\ 0\leq x \leq 1,\\
x+1 & \hbox{if }\ 1<x\leq 2,\\
1.5x & \hbox{if }\ 2<x.
\end{cases}
\end{equation*}
It is easy to verify that $h$ is continuous, strictly increasing and subadditive but not concave.

\begin{remark}\label{remark:3.3}
\emph{In Proposition \ref{prop3.1}, the condition that $h(u)\leq u$ holds for all $u\in[1,\infty]$ is essential when $s_{1}$ and $s_{2}$ are simultaneously normalized additive generators. Indeed, for any $x\in[1,\infty]$, let
\begin{equation*}
f_{1}(x)=k(x-1)+1 \mbox{ with } k>1,
\end{equation*}
\begin{equation*}
 f_{2}(x)=\begin{cases}
2x-1 & \hbox{if }\ 1\leq x \leq 2,\\
\frac{1}{2}x+2 & \hbox{if }\ 2<x\\
\end{cases}
\end{equation*}and
\begin{equation*}
f_{3}(x)=\ln(2e^{x}-e).
\end{equation*}
Then obviously, $f_{1}$, $f_{2}$ and $f_{3}$ are all continuous, strictly increasing, concave functions from $[1,\infty]$ to $[1,\infty]$ but none of them is subadditive.}
\end{remark}

Proposition \ref{prop3.1} enables us to determine the monotonicity of certain families of proper t-subnorms as illustrated by the following example.

\begin{figure}[!h]
\centering
\includegraphics[width=0.9\textwidth, keepaspectratio]{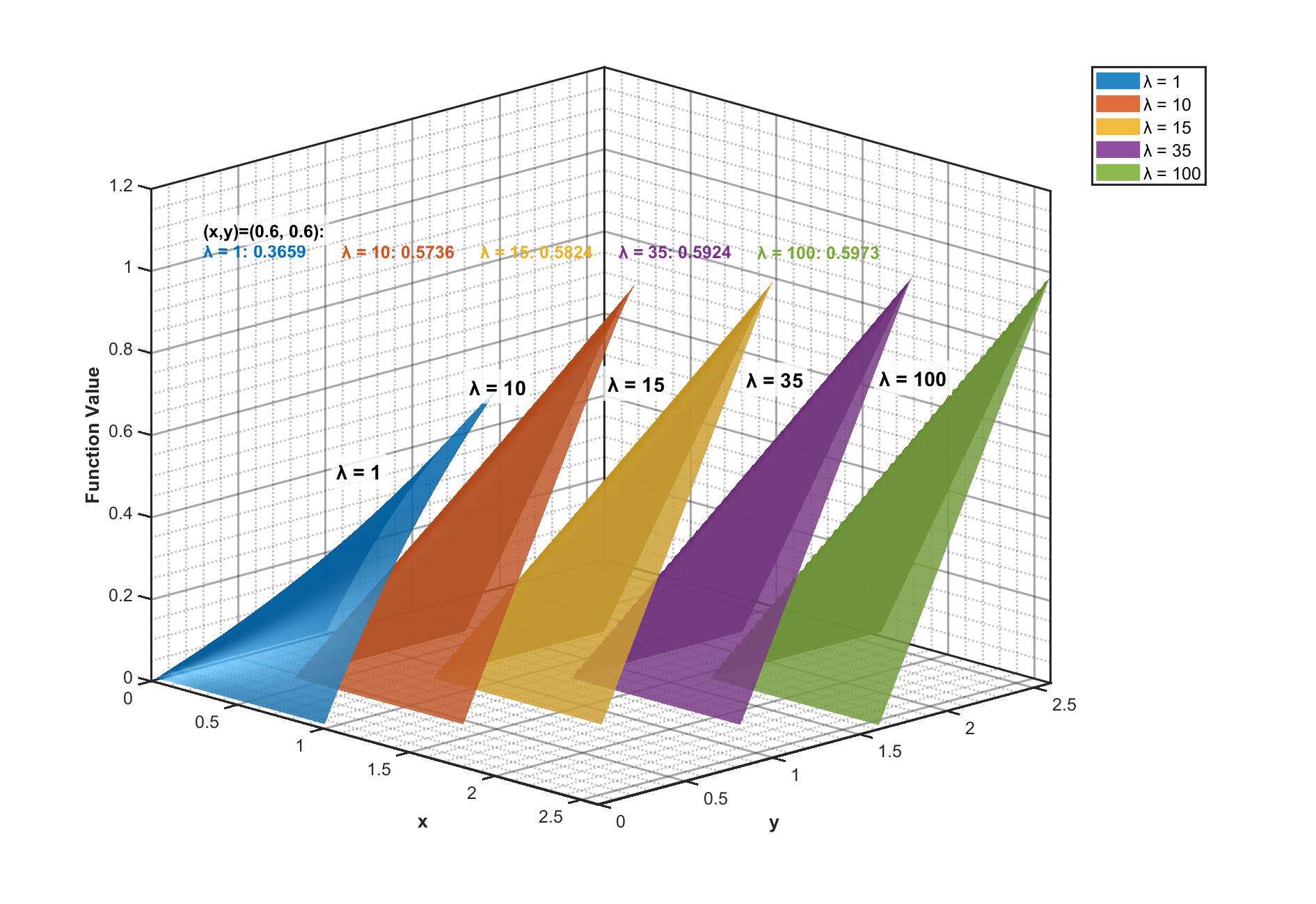}
\caption{A family of strictly increasing t-subnorms: $(S_{\lambda}^{\mathbf{D}})_{\lambda\in(0,\infty)}$ when $a=0.6$.}
\label{fig2}
\end{figure}

\begin{example}\label{exp3.3}\emph{
\renewcommand{\labelenumi}{(\roman{enumi})}
\begin{enumerate}
\item Consider the family $(T_{\lambda}^{\mathbf{D}})_{\lambda\in[0,\infty]}$ of \textit{Dombi $t$-norms} (see \cite{EP2000}). For any given $a\in (0,1)$, we can get a family of t-subnorms $(S_{\lambda}^{\mathbf{D}})_{\lambda\in (0,\infty)}\in\mathcal{S}_{p}$ defined by
    $$(S_{\lambda}^{\mathbf{D}})_{\lambda\in(0,\infty)}(x,y)=\frac{1}{a(1+((\frac{1-ax}{ax})^{\lambda}+(\frac{1-ay}{ay})^{\lambda})^{\frac{1}{\lambda}})}.$$
    Then their normalized additive generators are $s^\lambda:[0, 1]\rightarrow[1, \infty]$ with $s^{\lambda}(x)=(\frac{\frac{1}{x}-a}{1-a})^{\lambda}$, respectively.
    We claim that the family $(S_{\lambda}^{\mathbf{D}})_{\lambda\in(0,\infty)}$ is strictly increasing with respect to the parameter $\lambda$ (see Fig. \ref{fig2}). Indeed, for any  $\lambda, \mu \in (0, \infty)$ the function $s^\lambda \circ (s^\mu)^{-1}:[1, \infty]\longrightarrow[1, \infty]$ is given by $s^\lambda \circ (s^\mu)^{-1}(x)=x^{\frac{\lambda}{\mu}}$. When $\lambda \leq \mu$, this function is concave and satisfies $s^\lambda \circ (s^\mu)^{-1}(x)\leq x $ for all $x >1$. These properties together with Proposition \ref{prop3.1} show that $(S_{\lambda}^{\mathbf{D}})_{\lambda\in(0,\infty)}$ is a strictly increasing family of t-subnorms.
\item Similarly, for the family $(T_{\lambda}^{\mathbf{AA}})_{\lambda\in[0, \infty]}$ of Acz\'{e}l-Alsina t-norms, we derive a family of t-subnorms $(S_{\lambda}^{\mathbf{AA}})_{\lambda\in (0,\infty)}\in \mathcal{S}_{p}$ where $$(S_{\lambda}^{\mathbf{AA}})_{\lambda\in (0,\infty)}(x,y)=\frac{1}{a}\cdot e^{-((-\ln ax)^{\lambda} + (-\ln ay)^{\lambda})^{\frac{1}{\lambda}}}$$ with $a\in (0,1)$. This family is strictly increasing with respect to the parameter $\lambda $.
\end{enumerate}}
\end{example}

Below, let's recall that a function $f:I\rightarrow\mathbb{R}$, where $I\subseteq\mathbb{R}$ is an arbitrary non-empty interval, is called a \emph{convex function} if for all $x,y\in I$ and $\lambda\in[0,1]$ we have $f(\lambda x+(1-\lambda)y)\leq \lambda f(x)+(1-\lambda)f(y)$. Then we have the following proposition.

\begin{proposition}\label{prop3.2}
Let $S_{1}, S_{2}\in \mathcal{S}$ with additive generators $s_{1}:[0,1]\rightarrow[s_{1}(1), \infty]$ and $s_{2}:[0,1]\rightarrow[s_{2}(1), \infty]$, respectively. If the function $s_{1}\circ s_{2}^{-1}:[s_{2}(1), \infty]\rightarrow[s_{1}(1), \infty] $ is convex, then the following are equivalent:
\renewcommand{\labelenumi}{(\roman{enumi})}
\begin{enumerate}
\item $S_{1}\leq S_{2}$.
\item The function $s_{1}\circ s_{2}^{-1}$ is quasi-homogeneous, i.e., for all $x\in[s_{2}(1),\infty]$ and $t\geq 1$ we have
 $$s_{1}\circ s_{2}^{-1}(tx)\leq ts_{1}\circ s_{2}^{-1}(x).$$
\end{enumerate}
\end{proposition}
\begin{proof}
(i) $\Rightarrow$ (ii) If $S_{1}$ $\leq$ $S_{2}$, letting $h=s_{1}\circ s_{2}^{-1}$, then by Theorem \ref{theorem:3.1}, $h:[s_{2}(1), \infty]\rightarrow [s_{1}(1), \infty]$ is a continuous subadditive function. Letting $t\geq 1$, to prove that $h$ is quasi-homogeneous for all $x\in[s_{2}(1),\infty]$, we distinguish two cases as below.
\begin{itemize}
  \item If $t$ is a positive integer, then let $t=n$ with $n\in\mathbb{N}$. Obviously, the statement holds when $n=1$.  Now, assume that $h(kx)\leq kh(x)$ holds for all $x\in[s_{2}(1),\infty]$ when $n=k$. Then, when $n=k+1$, we have $h((k+1)x)=h(kx+x)\leq h(kx)+h(x)\leq kh(x)+h(x)=(k+1)h(x)$ since $h$ is subadditive. By induction, for any $n\in\mathbb{N}$ the statement $h(nx)\leq nh(x)$ holds. Therefore, $h$ is quasi-homogeneous.
  \item Next, consider the case that $t$ is not a positive integer. In this case, we can let $t=[t]+\{t\}$ where the notations $[t]$ and $\{t\}$ denote the integer and fractional parts of $t$, respectively. Subsequently, let $[t]=n$ with $n\in\mathbb{N}$. Then $t=n+\{t\}$ and $n\leq t <n+1$. Set $\lambda=n+1-t$. Then $0<\lambda\leq 1$ and $t=\lambda n+ (1-\lambda)(n+1)$. Thus we obtain
    \begin{align*} h(tx)&=h(\lambda nx+ (1-\lambda)(n+1)x)\\ &\leq \lambda h(nx)+(1-\lambda)h((n+1)x)\mbox{ from the convexity of }h\\ &\leq \lambda nh(x)+(1-\lambda)(n+1)h(x))\mbox{ since }h(nx)\leq nh(x)\mbox{ for all }x\in[s_{2}(1),\infty]\\ &=th(x).\end{align*}
    Therefore, $h$ is quasi-homogeneous.
\end{itemize}
(ii) $\Rightarrow$ (i) If $h$ is quasi-homogeneous, then by letting $t=2$, $h(2u)\leq 2h(u)$ for all $u\in[s_{2}(1),\infty]$. Because $h$ is a convex function, by setting $\lambda=\frac{1}{2}$ we obtain $h(\frac{u+v}{2})\leq \frac{1}{2}(h(u)+h(v))$. We distinguish three cases as follows.
\begin{itemize}
  \item If $u=v=s_{2}(1)$, then $h(u+v)=h(2s_{2}(1))\leq 2h(s_{2}(1))= 2s_{1}(1)=h(u)+h(v)$.
  \item If $u,v\in[s_{2}(1),\infty]$ with $s_{2}(1)<u+v<\infty$, then by the convexity and quasi-homogeneity of $h$ we have
  $$h(u+v)=h(2\cdot\frac{u+v}{2})\leq 2\cdot h(\frac{u+v}{2})\leq 2\cdot \frac{1}{2}(h(u)+h(v))=h(u)+h(v).$$
  \item If $u+v=\infty$, then
  $$h(u+v)=h(\infty)\leq h(\infty)+h(\min(u, v))=h(u)+h(v).$$
\end{itemize}
In summary, $h$ is subadditive, and $S_{1}\leq S_{2}$ by Theorem \ref{theorem:3.1}.
\end{proof}

Both Propositions \ref{prop3.1} and \ref{prop3.2} are applied to the comparison of continuous cancellative t-subnorms by virtue of the properties of the function $s_{1}\circ s_{2}^{-1}$. Next, we shall provide a more efficient criterion used for the comparison of continuous cancellative t-subnorms through the relationship between their respective additive generators $s_{1}$ and $s_{2}$.
\begin{proposition}\label{prop3.3}
 Let $S_{1}, S_{2}\in \mathcal{S}$ with additive generators $s_{1}:[0,1]\rightarrow[s_{1}(1), \infty]$ and $s_{2}:[0,1]\rightarrow[s_{2}(1), \infty]$, respectively. If $\frac{s_{1}}{s_{2}}$ is non-decreasing on $(0,1)$, then $S_{1}\leq S_{2}$.
\end{proposition}
\begin{proof}
Consider the function $\varphi:(s_{2}(1),\infty)\rightarrow(s_{1}(1),\infty)$ defined by
\begin{equation} \label{eq:1.12}
 \varphi(u)=\frac{s_{1}\circ s_{2}^{-1}(u)}{u}.
  \end{equation}
Fixing $x\in(0,1)$, and setting $u=s_{2}(x)$ in Eq.(\ref{eq:1.12}), we have $\varphi\circ s_{2}=\frac{s_{1}}{s_{2}}$. Since $\frac{s_{1}}{s_{2}}$ is non-decreasing and $s_{2}$ is strictly decreasing, $\varphi$ is non-increasing on $(s_{2}(1),\infty)$. Thus, for all $u,v\in(s_{2}(1),\infty)$
$$u[\varphi(u+v)-\varphi(u)]+v[\varphi(u+v)-\varphi(v)]\leq 0,$$
or equivalently,
$$(u+v)\varphi (u+v)\leq u\varphi(u)+v\varphi(v),$$
which together with Eq.(\ref{eq:1.12}) results in
$$s_{1}\circ s_{2}^{-1}(u+v)\leq s_{1}\circ s_{2}^{-1}(u)+s_{1}\circ s_{2}^{-1}(v).$$
Then it follows from the continuity of $s_{1}\circ s_{2}^{-1}$ that $s_{1}\circ s_{2}^{-1}$ is subadditive in $[s_{2}(1),\infty]$. Thus, $S_{1}\leq S_{2}$ by Theorem \ref{theorem:3.1}.
\end{proof}

\begin{example}\label{exp3.4}\emph{
\renewcommand{\labelenumi}{(\roman{enumi})}
\begin{enumerate}
\item Consider again $T_{1}$ and $T_{2}$ in Example \ref{exp3.2}$(\textup{ii})$. Since $\frac{t_{1}(x)}{t_{2}(x)}=-\frac{1}{\ln x}$ is increasing on $(0, 1)$, by Proposition \ref{prop3.3}, we have $T_{1}\leq T_{2}$.
\item For each $\lambda\in(0,\infty)$, the function $s^\lambda:[0,1]\rightarrow[0,\infty]$ is given by $s^{\lambda}(x)=(-\ln ax)^{\lambda}$ with $a\in(0,1]$. Obviously, $s^\lambda$ is an additive generator of a continuous cancellative t-subnorm denoted by $S^{(\lambda)}$. Moreover, the family $(S^{(\lambda)})_{\lambda\in(0,\infty)}$ is strictly increasing with respect to the parameter $\lambda$. Indeed, for $\lambda, \mu \in(0,\infty)$ with $\lambda< \mu$, we have
    $$\frac{s^{\lambda}(x)}{s^{\mu}(x)}=(-\ln ax)^{\lambda-\mu}.$$
    Let $g(x)=(-\ln ax)^{\lambda-\mu}$. Then $g$ is strictly increasing on $(0, 1)$, showing that $(S^{(\lambda)})_{\lambda\in(0,\infty)}$ is a strictly increasing family of t-subnorms by Proposition \ref{prop3.3}.
\end{enumerate}}
\end{example}

In general, the converse of Proposition \ref{prop3.3} is not true as shown by the following example.

\begin{example}\label{exp3.5}\emph{Consider the function $\psi:[1,\infty]\rightarrow[1,\infty]$ defined by
\begin{equation*}
 \psi(u)=\begin{cases}
-u^{2}+4u-2 & \hbox{if }\ 1\leq u\leq 2,\\
u & \hbox{if }\ u\geq 2.
\end{cases}
\end{equation*}
It is easy to see that $\psi$ is a strictly increasing bijection. Let $S_{2}$ $\in$ $\mathcal{S}_{p}$ with normalized additive generators $s_{2}$. Further, let $s_{1}=\psi\circ s_{2}$. Then the t-subnorm $S_{1}$ generated by $s_{1}$ also belongs to $\mathcal{S}_{p}$ since $\psi$ is a strictly increasing bijection. Moreover, the function $\psi$ is subadditive presented as follows.
\begin{itemize}
  \item If $1\leq u, v \leq 2$, then $2\leq u+v \leq 4$. In this case, $\psi(u)=-u^{2}+4u-2$, $\psi(v)=-v^{2}+4v-2$, and $\psi(u+v)=u+v$. Thus, we have
  $$\psi(u)+\psi(v)-\psi(u+v)=-u^{2}+3u-v^{2}+3v-4=-(u-\frac{3}{2})^{2}-(v-\frac{3}{2})^{2}+\frac{1}{2}\geq 0.$$
  \item If $1\leq u \leq 2$, $v\geq 2 $, then $u+v\geq 3$. In this case, $\psi(u)=-u^{2}+4u-2$, $\psi(v)=v$, and $\psi(u+v)=u+v$. Therefore,
  $$\psi(u)+\psi(v)-\psi(u+v)=-u^{2}+3u-2=-(u-2)(u-1)\geq 0.$$
 \item Similarly we can obtain $\psi(u)+\psi(v)-\psi(u+v)\geq 0$ when $1\leq v \leq 2$, $u\geq2$.
  \item If $u\geq 2$, $v\geq 2$, then $u+v\geq 4$. Thus $\psi(u)=u$, $\psi(v)=v$ and $\psi(u+v)=u+v$, i.e., $\psi(u)+\psi(v)=\psi(u+v)$.
\end{itemize}
 Meanwhile, it is easy to see that the subadditivity of $\psi$ is equivalent to that of $s_{1}\circ s_{2}^{-1}$. Thus, $S_{1}$ $\leq$ $S_{2}$ by Theorem \ref{theorem:3.1}. However, $\frac{s_{1}}{s_{2}}$ is not non-decreasing on $(0,1)$. Indeed, since $\frac{s_{1}}{s_{2}}=\frac{\psi\circ s_{2}}{s_{2}}$, we have
$\dfrac{s_1(x)}{s_2(x)}=\dfrac{\psi \circ s_2(x)}{s_2(x)}=- s_2(x)+4-\dfrac{2}{s_2(x)}$ for $ s_2(x)\in[1, 2]$. Then by a simple calculation, one can see that $\frac{s_{1}}{s_{2}}$ is strictly decreasing when $x\in[s_{2}^{-1}(\sqrt{2}), s_{2}^{-1}(1)]$.}
\end{example}

From the proof of Proposition \ref{prop3.3} the following corollary is immediately.

\begin{corollary}\label{coro3.3}
Let $S_{1}, S_{2}\in\mathcal{S}$ with additive generators $s_{1}:[0,1]\rightarrow[s_{1}(1), \infty]$ and $s_{2}:[0,1]\rightarrow[s_{2}(1), \infty]$, respectively. If the function $\varphi:(s_{2}(1), \infty)\rightarrow(s_{1}(1), \infty)$ defined by
\begin{equation*}
\label{eq:1.13}
\varphi(x)=\frac{s_{1}\circ s_{2}^{-1}(x)}{x}
\end{equation*}
is non-increasing, then $S_{1}\leq S_{2}$.
\end{corollary}

Analogously, the converse of Corollary \ref{coro3.3} does not hold generally.
%

It is very interesting that we have the following theorem.
\begin{theorem}\label{theorem:3.3}
Let $S_{1}, S_{2}\in\mathcal{S}$ with additive generators $s_{1}:[0,1]\rightarrow[s_{1}(1), \infty]$ and $s_{2}:[0,1]\rightarrow[s_{2}(1), \infty]$, respectively. If the function $s_{1}\circ s_{2}^{-1}:[s_{2}(1), \infty]\rightarrow[s_{1}(1), \infty]$ is convex, then the following are equivalent:
\renewcommand{\labelenumi}{(\roman{enumi})}
\begin{enumerate}
\item $S_{1}\leq S_{2}$.
\item The function $\frac{s_{1}\circ s_{2}^{-1}(x)}{x}$ is non-increasing on $(s_{2}(1), \infty)$.
\end{enumerate}
\end{theorem}
\begin{proof} Let $h=s_{1}\circ s_{2}^{-1}$. Then from Proposition \ref{prop3.2}, it suffices to prove that $h$ is quasi-homogeneous if and only if $\frac{h(x)}{x}$ is non-increasing.
If $h$ is quasi-homogeneous, then for all $x\in[s_{2}(1),\infty]$ and $t\geq 1$, we have $h(tx)\leq th(x)$. Let $0<x_{1}<x_{2}$, and $t=\frac{x_{2}}{x_{1}}$. Then
$$h(x_{2})=h(\frac{x_{2}}{x_{1}}\cdot x_{1})=h(t\cdot x_{1})\leq t\cdot h(x_{1})=\frac{x_{2}}{x_{1}}\cdot h(x_{1}),$$
or equivalently,
$$\frac{h(x_{2})}{x_{2}}\leq \frac{h(x_{1})}{x_{1}}.$$
Therefore, $\frac{h(x)}{x}$ is non-increasing.

Conversely, suppose that $\frac{h(x)}{x}$ is non-increasing, then for all $x\in(s_{2}(1), \infty)$ and $t\geq 1$, we have
$$\frac{h(tx)}{tx}\leq \frac{h(x)}{x}$$
since $tx\geq x$. Thus, $h(tx)\leq th(x)$, i.e., $h$ is quasi-homogeneous.
\end{proof}

The following proposition is related to the differentiability of the additive generators $s_1$ and $s_2$.
\begin{proposition}\label{prop3.4}
Let $S_{1}, S_{2}\in \mathcal{S}$ with differentiable additive generators $s_{1}:[0,1]\rightarrow[s_{1}(1), \infty]$ and $s_{2}:[0,1]\rightarrow[s_{2}(1), \infty]$, respectively. If the function $\frac{s_{1}'}{s_{2}'}$ is non-decreasing on $(0,1)$ and $s_{1}\leq s_{2}$ when both $s_{1}$ and $s_{2}$ are normalized additive generators, then $S_{1}\leq S_{2}$.
\end{proposition}
\begin{proof}
Let $h=s_{1}\circ s_{2}^{-1}$. Then
$$h'(u)=\frac{d}{du}[s_{1}\circ s_{2}^{-1}(u)]=s_{1}'(s_{2}^{-1}(u))\cdot \frac{d}{du}(s_{2}^{-1}(u)).$$
Let $u=s_{2}(x)$. Then we have
$$\frac{d}{du}(s_{2}^{-1}(u))=\frac{1}{s'_{2}(x)}=\frac{1}{s'_{2}(s_{2}^{-1}(u))}$$
since $s_{2}$ is continuous differentiable and strictly decreasing. Hence,
$$h'(u)=\frac{s_{1}'}{s_{2}'}(s_{2}^{-1}(u)).$$
Moreover, as the function $\frac{s_{1}'}{s_{2}'}$ is non-decreasing on $(0,1)$ and $s_{2}^{-1}$ is strictly decreasing in $[s_{2}(1),\infty]$, we have that the function $h'$ is non-increasing. Consequently, the function $h=s_{1}\circ s_{2}^{-1}:[s_{2}(1), \infty]\rightarrow[s_{1}(1), \infty]$ is concave. Since $u=s_{2}(x)$ and $s_{1}(x)\leq s_{2}(x)$ for all $x\in[0,1]$ when both $s_{1}$ and $s_{2}$ are normalized additive generators, we have $s_{1}\circ s_{2}^{-1}(u)\leq u$ for all $u\in[1,\infty]$. It follows from Proposition \ref{prop3.1} that $S_{1}\leq S_{2}$.
\end{proof}

Proposition \ref{prop3.4} can also be used to determine whether a families of proper t-subnorms is monotone or not.
\begin{example}\label{exp3.7}\emph{For the family $(T_{\lambda}^{\mathbf{SS}})_{\lambda\in[-\infty, \infty]}$ of Schweizer-Sklar t-norms (see \cite{EP2000}), when $\lambda\in(-\infty,0]$, a family of t-subnorms $(S_{\lambda}^{\mathbf{SS}})_{\lambda\in (-\infty,0]}\in \mathcal{S}_{p}$ defined by
    $$(S_{\lambda}^{\mathbf{SS}})_{\lambda\in(-\infty,0]}(x,y)=\frac{1}{a}(\max((ax)^{\lambda}+(ay)^{\lambda}-1, 0))^{\frac{1}{\lambda}}$$
    can be obtained where $a\in(0,1)$. Then their normalized additive generators are $s:[0,1]\rightarrow[1,\infty]$ with $s^{\lambda}(x)=\frac{1-(ax)^{\lambda}}{1-a^{\lambda}}$, respectively.
    Obviously, for each $x\in(0,1)$ and all $\lambda,\mu\in(-\infty,0]$ with $\lambda< \mu$, we always have
    $$\frac{s_{\mu}'(x)}{s_{\lambda}'(x)}=\frac{\mu(1-a^{\lambda})}{\lambda(1-a^{\mu})}\cdot (ax)^{\mu-\lambda},$$
    which is non-decreasing on $(0,1)$. Moreover, simple calculations show that $s_{\mu}\leq s_{\lambda}$. Then $(S_{\mu}^{\mathbf{SS}})\leq (S_{\lambda}^{\mathbf{SS}})$ from Proposition \ref{prop3.4}, i.e., $(S_{\lambda}^{\mathbf{SS}})_{\lambda\in(-\infty,0]}$ is a strictly decreasing family of t-subnorms.}
\end{example}

It is well-known that each t-subnorm $S$ can be transformed into a t-norm by redefining (if necessary) its values on the upper right boundary of the unit square as the following proposition.

\begin{proposition}[\cite{EP2000}]\label{prop3.5}
If $S$ is a t-subnorm then the function $T:[0,1]^2\rightarrow[0,1]$ defined by
\begin{equation*}
T(x,y) =
\begin{cases}
S(x,y) & \text{if } (x,y) \in [0,1)^2, \\
\min(x,y) & \text{otherwise},
\end{cases}
\end{equation*}
is a t-norm.
\end{proposition}

Obviously, $S$ and $T$ in Proposition \ref{prop3.5} satisfy $S\leq T$. Further, let $T_{1}$ be the t-norm generated by the t-subnorm $S$ via Proposition \ref{prop3.5}, and let $T_{2}$ be an arbitrary t-norm. Then $T_{1}\leq T_{2}$ if and only if $S\leq T_{2}$.

 In what follows, we investigate the order relation between a proper t-subnorm and a t-norm. Specifically, we shall restrict ourselves to the case when $S\in\mathcal{S}_{p}$ and $T$ is a continuous Archimedean t-norm. We first recall a result from \cite{EP2000}.
\begin{lemma}[\cite{EP2000}]\label{lemm3.1}
For a function $T:[0,1]^2\rightarrow[0,1]$ the following are equivalent:
\renewcommand{\labelenumi}{(\roman{enumi})}
\begin{enumerate}
    \item $T$ is a strict  t-norm.
    \item $T$ is isomorphic to the product $T_{\mathbf{P}}$, i.e., there is a strictly increasing bijection $\varphi:[0,1]\rightarrow[0,1]$ such that for all $(x,y) \in [0,1]^2 $
  $$T(x,y)=\varphi^{-1}(T_{\mathbf{P}}(\varphi(x), \varphi(y)))=\varphi^{-1}(\varphi(x)\cdot \varphi(y)).$$
\end{enumerate}
\end{lemma}

Then we have the following theorem that is different from Remark \ref{remark:3.1}$(\textup{ii})$.
\begin{theorem}\label{theorem:3.4}
Let $S\in\mathcal{S}_{p}$ with the normalized additive generator $s:[0,1]\rightarrow[1, \infty]$, and $T\in\mathcal{T}_{s}$. Then the following are equivalent:
\renewcommand{\labelenumi}{(\roman{enumi})}
\begin{enumerate}
\item $S\leq T$.
\item There is a strictly increasing bijection $\varphi:[0,1]\rightarrow[0,1]$ such that the function $s\circ \varphi^{-1}$ is submultiplicative-additive, i.e., for all $u,v\in[0,1]$ we have
    $$s\circ \varphi^{-1}(uv)\leq s\circ \varphi^{-1}(u)+s\circ \varphi^{-1}(v).$$
\end{enumerate}
\end{theorem}
\begin{proof}
Since $S\in \mathcal{S}_{p}$, by Eq.(\ref{eq:1.1}) we have
\begin{equation} \label{eq:1.14}
S(x,y)=s^{-1}(s(x)+s(y)).
\end{equation}
Moreover, as a consequence of Lemma \ref{lemm3.1}, there is a strictly increasing bijection $\varphi:[0,1]\rightarrow[0,1]$ such that for all $(x,y)\in[0,1]^2$
\begin{equation} \label{eq:1.15}
 T(x,y)=\varphi^{-1}(T_{\mathbf{P}}(\varphi(x), \varphi(y)))=\varphi^{-1}(\varphi(x) \cdot \varphi(y))
\end{equation}
and $t(x)=t_{p}\circ \varphi(x)=-\ln \varphi(x)$ is an additive generator of $T$.
Thus, from Eqs.(\ref{eq:1.14}) and (\ref{eq:1.15}) we can obtain $S\leq T$ if and only if for all $x,y\in[0,1]$
\begin{equation} \label{eq:1.16}
s^{-1}(s(x)+s(y))\leq \varphi^{-1}(\varphi(x)\cdot \varphi(y)).
\end{equation}
Fixing $x,y\in[0,1]$, setting $u=\varphi(x)$ and $v=\varphi(y)$, and applying the strictly decreasing bijection $s$ to both sides of Eq.(\ref{eq:1.16}) we have
$$s\circ \varphi^{-1}(uv)\leq s\circ \varphi^{-1}(u)+s\circ \varphi^{-1}(v).$$
Thus, $S\leq T$ if and only if $s\circ \varphi^{-1}$ is submultiplicative-additive.
\end{proof}

Theorem \ref{theorem:3.4} enables the construction of a strict t-norm that is stronger than a given continuous cancellative proper t-subnorm, and the strict t-norm constructed by this method differs from that obtained via Proposition \ref{prop3.5}.

\begin{example}\label{exp3.8}\emph{Consider the function $S:[0, 1]^2\rightarrow[0, 1]$ with $S(x,y)=\frac{1}{2}xy$. Then $S\in \mathcal{S}_{p}$ and its normalized additive generator of $S$ is the function $s:[0, 1]\rightarrow[1, \infty]$ with $s(x)=1-\frac{\ln x}{\ln 2}$. For the strictly increasing bijection $\varphi:[0,1]\rightarrow[0,1]$ defined by $\varphi(x)=e^{1-\frac{1}{x}}$, we obtain that $ s\circ \varphi^{-1}(x)=1+\frac{\ln (1-\ln x)}{\ln 2}$.
A simple computation shows that
$$ s\circ \varphi^{-1}(xy)\leq s\circ \varphi^{-1}(x)+s\circ \varphi^{-1}(y)$$
for all $x,y\in[0,1]$. From Theorem \ref{theorem:3.4}, there exists a strict t-norm $T$ such that $S\leq T$ that is exactly the Hamacher product $T_{0}^{H}$ with the additive generator $t(x)=\frac{1-x}{x}$.}
\end{example}

Furthermore, we have the following theorem.
\begin{theorem}\label{theorem:3.5}
Let $S\in \mathcal{S}$ with an additive generator $s:[0,1]\rightarrow[s(1), \infty]$, and $T\in \mathcal{T}_{s}$. Then the following are equivalent:
\renewcommand{\labelenumi}{(\roman{enumi})}
\begin{enumerate}
\item $S=T$.
\item There is a strictly increasing bijection $\varphi: [0,1]\rightarrow[0,1]$ such that the function $s\circ \varphi^{-1}$ is logarithmic.
\end{enumerate}
\end{theorem}
\begin{proof}
From the proof of Theorem \ref{theorem:3.4}, we know that $S=T$ if and only if for all $x,y\in[0,1]$
$$ s^{-1}(s(x)+s(y))=\varphi^{-1}(\varphi(x) \cdot \varphi(y)),$$
or equivalently,
\begin{equation} \label{eq:1.17}
 s\circ \varphi^{-1}(uv)= s\circ \varphi^{-1}(u)+s\circ \varphi^{-1}(v),
\end{equation}
where $u,v\in[0,1]$ with $u=\varphi(x)$ and $v=\varphi(y)$. Evidently, Eq.(\ref{eq:1.17}) is a Cauchy-type equation, whose continuous, strictly decreasing solutions $ s\circ \varphi^{-1}:[0,1]\rightarrow[0,\infty]$ must satisfy $s\circ \varphi^{-1}(x)=-c\cdot \ln x$ for some $c\in(0, \infty)$. Therefore, $S=T$ if and only if $s\circ \varphi^{-1}(x)=-c\cdot \ln x$ for some $c\in(0, \infty)$.
\end{proof}

As a conclusion of this section we consider the case that $T$ is a nilpotent t-norm.
\begin{proposition}\label{prop3.6}
Let $S\in\mathcal{S}$. Then there exists no nilpotent t-norm $T$ such that $S\leq T$.
\end{proposition}
\begin{proof}
 Let $T$ be a nilpotent t-norm. Then for each $x\in(0,1)$, there exists an $n\in \mathbb{N}$ such that $x_{T}^{(n)}=0$. If $S\leq T$, then $x_{S}^{(n)}\leq x_{T}^{(n)}=0$, implying $x_{S}^{(n)}=0$, contrary to the fact that $S\in\mathcal{S}$.
\end{proof}

\section{Comparison of continuous cancellative t-subnorms based on growth and boundedness}
This section explores the growth and boundedness of additive generators when two continuous cancellative t-subnorms are comparable, which are applied to the comparison of continuous cancellative t-subnorms.

We first cite the following lemma from \cite{HP1957}.
\begin{lemma}[\cite{HP1957}]\label{lemm4.1}
If $f(x)$ is subadditive and finite on $(a,\infty)$, $a\geq 0$, then
$$ \lim_{x\rightarrow\infty}\frac{f(x)}{x}=\inf_{x>a}\frac{f(x)}{x}<\infty.$$
\end{lemma}

\begin{proposition}\label{prop4.1}
Let $S_{1}, S_{2}\in\mathcal{S}$ with additive generators $s_{1}:[0,1]\rightarrow[s_{1}(1), \infty]$ and $s_{2}:[0,1]\rightarrow[s_{2}(1), \infty]$, respectively.  If $S_{1}\leq S_{2}$, then
$$\lim_{x\rightarrow\infty}\frac{s_{1}\circ s_{2}^{-1}(x)}{x}=\inf_{x>s_{2}(1)}\frac{s_{1}\circ s_{2}^{-1}(x)}{x}<\infty.$$
\end{proposition}
\begin{proof}
Since $S_{1}\leq S_{2}$, by Theorem \ref{theorem:3.1} we know that $s_{1}\circ s_{2}^{-1}$ is continuous, strictly increasing and subadditive. It follows from Lemma \ref{lemm4.1} that
$$\lim_{x\rightarrow\infty}\frac{s_{1}\circ s_{2}^{-1}(x)}{x}=\inf_{x>s_{2}(1)}\frac{s_{1}\circ s_{2}^{-1}(x)}{x}<\infty.$$
\end{proof}

Proposition \ref{prop4.1} shows that for two continuous cancellative t-subnorms $S_{1}$ and $S_{2}$ with $S_{1}\leq S_{2}$, the average growth rate $\frac{s_{1}\circ s_{2}^{-1}(x)}{x}$ of $s_{1}\circ s_{2}^{-1}$ converges to a finite constant. This captures the asymptotic linearity of $s_{1}\circ s_{2}^{-1}$, i.e., $s_{1}\circ s_{2}^{-1}$ is dominated above by a suitably chosen linear function of $x$ for large positive values of $x$. Moreover, $\displaystyle\lim_{x\rightarrow\infty}\frac{s_{1}\circ s_{2}^{-1}(x)}{x}=\displaystyle\lim_{t\rightarrow 0^{+}}\frac{s_{1}(t)}{s_{2}(t)}$. This implies the following two statements, which can be used for the construction of two comparable continuous cancellative t-subnorms:
\begin{itemize}
\item If $\displaystyle\lim_{t\rightarrow 0^{+}}\frac{s_{1}(t)}{s_{2}(t)}=0$, then $s_{1}$ is an infinity of lower order than $s_{2}$ as $t\rightarrow 0^{+}$.
\item If $0<\displaystyle\lim_{t\rightarrow 0^{+}}\frac{s_{1}(t)}{s_{2}(t)}<\infty$, then $s_{1}$ and $s_{2}$ are infinities of same order as $t\rightarrow 0^{+}$.
\end{itemize}

Then we have the following theorem.
\begin{theorem}\label{theorem:4.1}
 Let $S_{1}, S_{2}\in\mathcal{S}$ with additive generators $s_{1}:[0,1]\rightarrow[s_{1}(1), \infty]$ and $s_{2}:[0,1]\rightarrow[s_{2}(1), \infty]$, respectively. If the function $\frac{s_{1}\circ s_{2}^{-1}(x)}{x}$ is convex on $(s_{2}(1),\infty)$, then the following are equivalent:
\renewcommand{\labelenumi}{(\roman{enumi})}
\begin{enumerate}
\item $S_{1}\leq S_{2}$.
\item $\displaystyle\lim_{x\rightarrow\infty}\frac{s_{1}\circ s_{2}^{-1}(x)}{x}=\displaystyle\inf_{x>s_{2}(1)}\frac{s_{1}\circ s_{2}^{-1}(x)}{x}<\infty$.
\end{enumerate}
\end{theorem}
\begin{proof}
(i) $\Rightarrow$ (ii) The proof is immediate from Proposition \ref{prop4.1}.

(ii) $\Rightarrow$ (i) Let $\varphi(x)=\frac{s_{1}\circ s_{2}^{-1}(x)}{x}$. Then by Corollary \ref{coro3.3}, it suffices to prove that $\varphi(x)$ is non-increasing. Assume there exists $a,b\in(s_{2}(1),\infty)$ with $a<b$ such that $\varphi(a)<\varphi(b)$. For any $x>b$, let $\lambda=\frac{x-b}{x-a}$. Then $b=\lambda a+(1-\lambda)x$. Thus
$$ \varphi(b)=\varphi(\lambda a+(1-\lambda)x)\leq \lambda\varphi(a)+(1-\lambda)\varphi(x)$$
since $\varphi$ is convex. Hence
\begin{equation*}
 \varphi(x)\geq\frac{x(\varphi(b)-\varphi(a))+b\varphi(a)-a\varphi(b)}{b-a}.
\end{equation*}
 This follows that $\displaystyle\lim_{x\rightarrow\infty}\varphi(x)=\infty$ since $\varphi(a)<\varphi(b)$, which contradicts $(\textup{ii})$.
\end{proof}

\begin{lemma} [\cite{MK2009}]\label{lemm4.2}
Let $f:\mathbb{R}\rightarrow\mathbb{R}$ be a measurable subadditive function, and let $A=\displaystyle\inf_{x<0}\frac{f(x)}{x}$ (resp. $B=\displaystyle\sup_{x>0}\frac{f(x)}{x}$). If $A$ (resp. $B$) is finite, then
$$ A=\lim_{x\rightarrow0^{-}}\frac{f(x)}{x} \mbox{ (resp. } B=\lim_{x\rightarrow0^{+}}\frac{f(x)}{x}).$$
\end{lemma}

From Lemma \ref{lemm4.2} and Theorem \ref{theorem:3.1}, we immediately have the following proposition.
\begin{proposition}\label{prop4.4}
 Let $T_{1}, T_{2}\in\mathcal{T}_{s}$ with additive generators $t_{1}, t_{2}:[0,1]\rightarrow[0, \infty]$, respectively. If $T_{1}\leq T_{2}$, then $\frac{t_{1}\circ t_{2}^{-1}(x)}{x}$ is bounded on $(0,\infty)$ if and only if $Ax\leq t_{1}\circ t_{2}^{-1}(x)\leq Bx$, where $A=\displaystyle\lim_{x\rightarrow\infty}\frac{t_{1}\circ t_{2}^{-1}(x)}{x}$ and $B=\displaystyle\lim_{x\rightarrow 0^+}\frac{t_{1}\circ t_{2}^{-1}(x)}{x}$.
\end{proposition}

\begin{corollary}\label{coro4.2}
Let $S\in\mathcal{S}_{p}$, $T\in\mathcal{T}_{s}$ with additive generators $s:[0,1]\rightarrow[1, \infty]$ and $t:[0,1]\rightarrow[0, \infty]$, respectively. If $S\leq T$, then $\frac{s\circ t^{-1}(x)}{x}$ is bounded on $(0,\infty)$ if and only if $Ax\leq s\circ t^{-1}(x)\leq Bx$, where $A=\displaystyle\lim_{x\rightarrow\infty}\frac{s\circ t^{-1}(x)}{x}$ and $B=\displaystyle\lim_{x\rightarrow 0^+}\frac{s\circ t^{-1}(x)}{x}$.
\end{corollary}

\begin{remark}\label{remark:4.1}
\emph{In Proposition \ref{prop4.4} and Corollary \ref{coro4.2}, if $\frac{h(x)}{x}$ has no upper bounds, then $h$ has also no upper bounds.} 
\end{remark}

\begin{proposition}\label{prop4.5}
Let $S_{1}, S_{2}\in\mathcal{S}_{p}$ with normalized additive generators $s_{1}, s_{2}:[0,1]\rightarrow[1, \infty]$, respectively. If $S_{1}\leq S_{2}$, then
$Ax\leq s_{1}\circ s_{2}^{-1}(x)\leq Bx$, where $A=\displaystyle\lim_{x\rightarrow\infty}\frac{s_{1}\circ s_{2}^{-1}(x)}{x}$ and $B=\displaystyle\sup_{x\geq1}\frac{s_{1}\circ s_{2}^{-1}(x)}{x}$.
\end{proposition}
\begin{proof}
Let $h=s_{1}\circ s_{2}^{-1}$ and $H=\{\frac{h(x)}{x}| x\geq 1\}$. Then $H$ is non-empty since $h(1)\in H$. Thus by Proposition \ref{prop4.1}, $A=\displaystyle\lim_{x\rightarrow \infty}\frac{h(x)}{x}=\inf H$. Next, we show that $H$ has an upper bound by distinguishing two cases as below.
\begin{itemize}
  \item If $x=n\in\mathbb{N}$, then $h(n)\leq nh(1)=n$ since $h$ is subadditive and $h(1)=1$, i.e., $\frac{h(n)}{n}\leq1$.
  \item If $n< x<n+1$ with $n\in \mathbb{N}$, then $\frac{h(x)}{x}< \frac{n+1}{n}\leq2$ for all $n\geq1$ since $h$ is strictly increasing and $h(n)\leq n$.
\end{itemize}
Thus $2$ is an upper bound of $H$. Therefore, $B=\sup H$ exists.
\end{proof}

From Proposition \ref{prop4.4} together with Corollaries \ref{coro3.3} and \ref{coro4.2} we easily deduce the following theorem.
\begin{theorem}\label{theorem:4.2}
Let $S_{1}\in\mathcal{S}$, $S_{2}\in\mathcal{T}_{s}$ with additive generators $s_{1}:[0,1]\rightarrow[s_{1}(1), \infty]$ and $s_{2}:[0,1]\rightarrow[0, \infty]$, respectively. If the function $\frac{s_{1}\circ s_{2}^{-1}(x)}{x}$ is non-increasing and bounded on $(0,\infty)$, then
\renewcommand{\labelenumi}{(\roman{enumi})}
\begin{enumerate}
\item $S_{1}\leq S_{2}$.
\item $Ax\leq s_{1}\circ s_{2}^{-1}(x)\leq Bx$, where $A=\displaystyle\lim_{x\rightarrow\infty}\frac{s_{1}\circ s_{2}^{-1}(x)}{x}$ and $B=\displaystyle\lim_{x\rightarrow 0^{+}}\frac{s_{1}\circ s_{2}^{-1}(x)}{x}$.
\end{enumerate}
\end{theorem}

Furthermore, we have the following theorem.
\begin{theorem}\label{theorem:4.3}
Let $S_{1}, S_{2}\in\mathcal{S}$ with additive generators $s_{1}:[0,1]\rightarrow[s_{1}(1), \infty]$ and $s_{2}:[0,1]\rightarrow[s_{2}(1), \infty]$, respectively. If the function $s_{1}\circ s_{2}^{-1}:[s_{2}(1), \infty]\rightarrow[s_{1}(1), \infty]$ is concave and $\displaystyle\sup_{x>s_{2}(1)}\frac{s_{1}\circ s_{2}^{-1}(x)}{x}\leq 1$, then the following are equivalent:
\renewcommand{\labelenumi}{(\roman{enumi})}
\begin{enumerate}
\item $S_{1}\leq S_{2}$.
\item $A\leq \frac{s_{1}\circ s_{2}^{-1}(x)}{x}\leq 1$, where $A=\displaystyle\lim_{x\rightarrow\infty}\frac{s_{1}\circ s_{2}^{-1}(x)}{x}$.
\end{enumerate}
\end{theorem}
\begin{proof}
(i) $\Rightarrow$ (ii) If $S_{1}\leq S_{2}$, then by Proposition \ref{prop4.1} we have $s_{1}\circ s_{2}^{-1}(x)\geq Ax$, where $A=\displaystyle\lim_{x\rightarrow\infty}\frac{s_{1}\circ s_{2}^{-1}(x)}{x}=\displaystyle\inf_{x>s_{2}(1)}\frac{s_{1}\circ s_{2}^{-1}(x)}{x}$. On the other hand, $s_{1}\circ s_{2}^{-1}(x)\leq x$ for all $x\in[s_{2}(1), \infty]$ since $\displaystyle\sup_{x>s_{2}(1)}\frac{s_{1}\circ s_{2}^{-1}(x)}{x}\leq1$. Thus $Ax\leq s_{1}\circ s_{2}^{-1}(x)\leq x$, i.e., $A\leq \frac{s_{1}\circ s_{2}^{-1}(x)}{x}\leq 1$.

(ii) $\Rightarrow$ (i) If $A\leq \frac{s_{1}\circ s_{2}^{-1}(x)}{x} \leq 1$, then $Ax\leq s_{1}\circ s_{2}^{-1}(x)\leq x$. Since $h$ is concave, we get $S_{1}\leq S_{2}$ by Proposition \ref{prop3.1}.
\end{proof}

\begin{figure}[!h]
    \centering
    \begin{subfigure}[t]{0.9\textwidth}
        \centering
        \includegraphics[width=\textwidth]{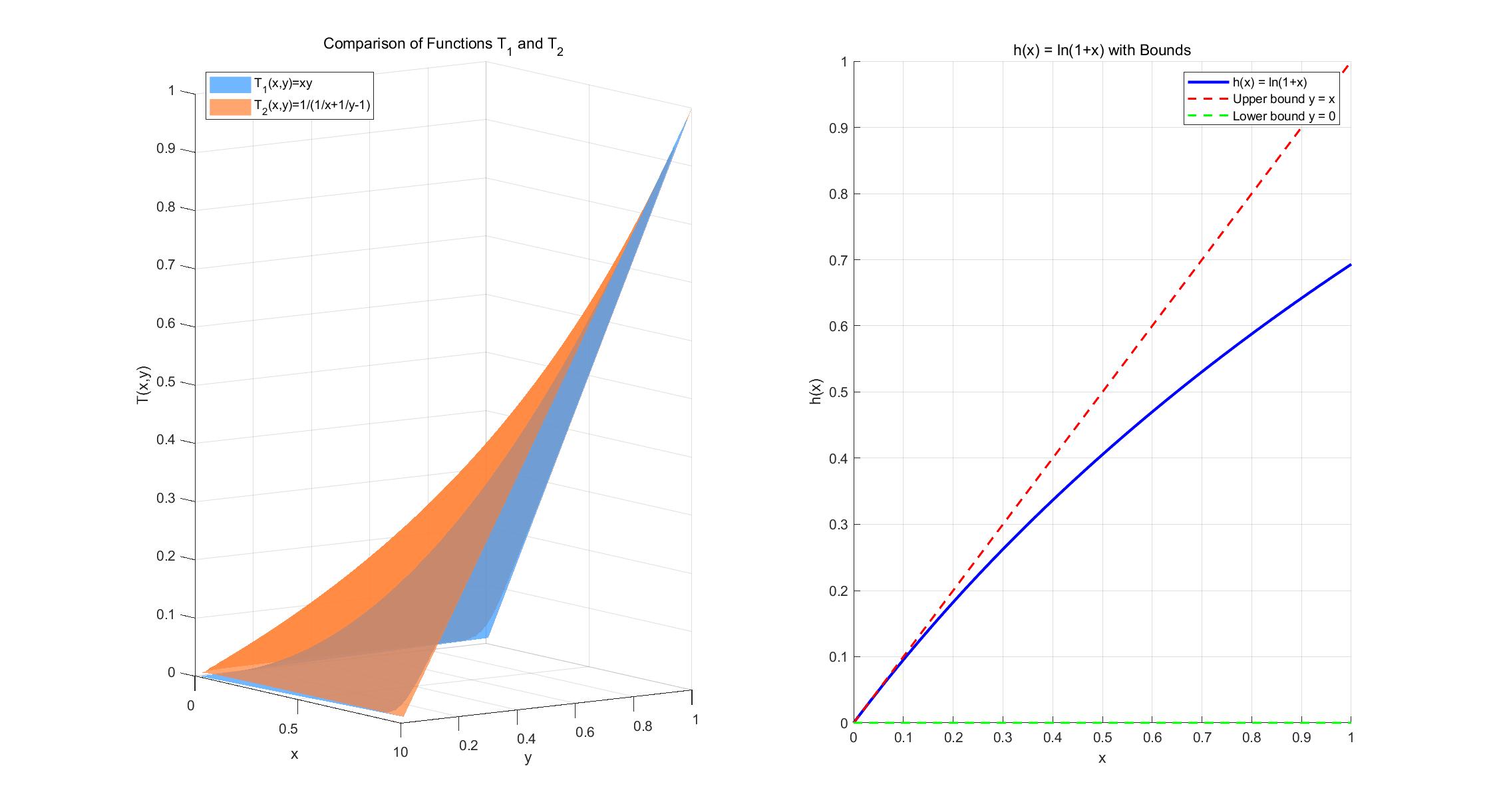}
        \caption{$T_{1} \leq T_{2}$ (see Example \ref{exp4.1}$(\textup{i})$).}
    \end{subfigure}
    \hfill
    \begin{subfigure}[t]{0.9\textwidth}
        \centering
        \includegraphics[width=\textwidth]{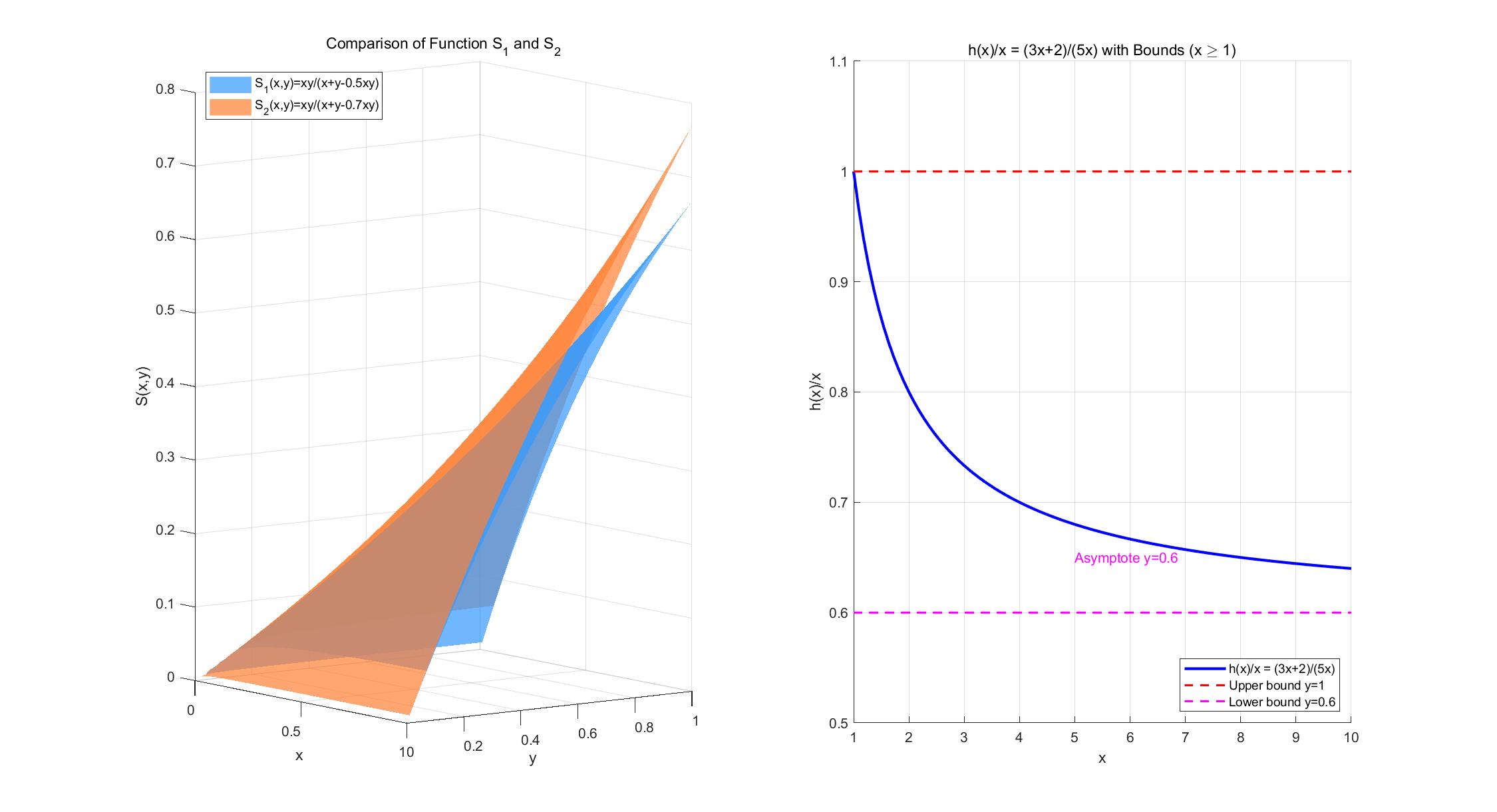}
        \caption{$S_{1} \leq S_{2}$ (see Example \ref{exp4.1}$(\textup{ii})$).}
    \end{subfigure}
    \caption{3D plots of comparable t-subnorms ((a) $T_{1}, T_{2}\in \mathcal{T}_{s},$; (b) $S_{1},S_{2}\in\mathcal{S}_{p}$).}
    \label{fig3}
\end{figure}

The following example illustrates Theorems \ref{theorem:4.2} and \ref{theorem:4.3}.

\begin{example}\label{exp4.1}\emph{
\renewcommand{\labelenumi}{(\roman{enumi})}
\begin{enumerate}
\item In Example \ref{exp3.1}(i), $T_{1}=T_{p}, T_{2}=T_{0}^{H}$ and $t_{1}\circ t_{2}^{-1}(x)=\ln(x+1)$. Obviously, the function $\frac{\ln(1+x)}{x}$ is decreasing, $A=\displaystyle\lim_{x\rightarrow\infty}\frac{\ln(x+1)}{x}=0$ and $B=\displaystyle\lim_{x\rightarrow 0^{+}}\frac{\ln(x+1)}{x}=1$. Hence $T_{1}\leq T_{2}$ and $0\leq t_{1}\circ t_{2}^{-1}(x)\leq x$ by Theorem \ref{theorem:4.2} (see Fig. \ref{fig3}(a)).
\item In Example \ref{exp3.1}(iii), let $h(x)=s_{1}\circ s_{2}^{-1}(x)$. Then $h(x)=\frac{3x+2}{5}$. Evidently, $h$ is concave and the function $\frac{}{}\frac{}{}\frac{h(x)}{x}=\frac{3x+2}{5x}$ is strictly decreasing. Then $\displaystyle\sup_{x\geq1}\frac{3x+2}{5x}=\displaystyle\lim_{x\rightarrow 1^{+}}\frac{3x+2}{5x}=1$. Moreover, $A=\displaystyle\lim_{x\rightarrow\infty}\frac{3x+2}{5x}=\frac{3}{5}$. So that $\frac{3}{5}\leq \frac{h(x)}{x} \leq 1$. Thus $S_{1}\leq S_{2}$ by Theorem \ref{theorem:4.3} (see Fig. \ref{fig3}(b)).
\end{enumerate}}
\end{example}

\section{Conclusions}
 This article focuses on comparing continuous cancellative t-subnorms. For continuous cancellative t-subnorms, it not only established some necessary and sufficient conditions for their comparability but also supplied several practical sufficient criteria, which are applied to determine the monotonicity of some families of t-subnorms. In addition, this article related growth and boundedness of additive generators of continuous cancellative t-subnorms to their order relations.


\end{document}